\documentclass [11pt,oneside]{article}

\reversemarginpar \textwidth=14cm \textheight=19cm

\newcommand{\figfile}[1]{{\Large #1}}
\newcommand{\figfatta}[1]{\includegraphics{#1}}

\usepackage{amssymb}
\usepackage{amsfonts}
\usepackage{amsmath}
\usepackage{amsthm}
\usepackage[dvips]{graphicx}

\newtheorem{lemma}{Lemma}[section] 
\newtheorem{teo}[lemma]{Theorem}
\newtheorem{rem}[lemma]{Remark} 
\newtheorem{prop}[lemma]{Proposition}
\newtheorem{cor}[lemma]{Corollary}

\newcommand{\matN}{\ensuremath {\mathbb{N}}}

\newcommand{\matZ} {\ensuremath {\mathbb{Z}}}

\newcommand{\matP} {\ensuremath {\mathbb{P}}}

\newcommand{\calX} {\ensuremath {\mathcal{X}}}

\newcommand{\calL} {\ensuremath {\mathcal{L}}}

\newcommand{\calM} {\ensuremath {\mathcal{M}}}

\newcommand{\calB} {\ensuremath {\mathcal{B}}}

\newcommand{\calF} {\ensuremath {\mathcal{F}}}

\newcommand{\nota} [1] {\caption{\footnotesize{#1}}}

\def\interior#1{{\rm int}(#1)}

\def\closedpair#1{#1}

\font\titsc=cmcsc10 scaled 1200

\newcommand{\ptwoirred}{$\matP^2$-irreducible}
\newcommand{\solK}{\hbox{\textit{\textbf{K}}}}
\newcommand{\solT}{\hbox{\textit{\textbf{T}}}}

\newcommand{\timtil}{\begin{picture}(12,12)
\put(2,0){$\times$}\put(2,4.5){$\sim$}\end{picture}}

\newcommand{\finedimo}{{\hfill\hbox{$\square$}\vspace{2pt}}}
\newcommand{\dimo}[1]{\vspace{2pt}\noindent\textit{Proof of #1}.}

\author{Bruno \titsc{Martelli} \and Carlo \titsc{Petronio}}

\title{A New Decomposition Theorem for 3-Manifolds}

\begin{document}

\maketitle

\tableofcontents

\section*{Introduction}
We develop in this paper a theory of complexity for pairs $(M,X)$ where $M$ is
a compact 3-manifold such that $\chi(M)=0$,
and $X$ is a collection of trivalent graphs, each graph $\tau$ being embedded in one
component $C$ of $\partial M$ so that $C\setminus\tau$ is one disc.
In the special case where $M$ is closed, so $X=\emptyset$, our complexity coincides with
Matveev's~\cite{Matveev:AAM}. Extending his results
we show that complexity of pairs is additive under connected sum and that,
when $M$ is closed, irreducible, \ptwoirred\
and different from $S^3,L_{3,1},\matP^3$, its complexity
is precisely the minimal number of tetrahedra in a triangulation.
These two facts show that indeed complexity is a very natural measure
of how complicated a manifold or pair is.
The former fact was known to Matveev in the closed case,
the latter one in the orientable case.

The most relevant feature of our
theory is that it leads to a splitting theorem along tori and Klein bottles
for irreducible and \ptwoirred\ pairs (so, in particular, for irreducible and \ptwoirred\
closed manifolds). The blocks of the splitting are themselves pairs, and the complexity
of the original pair is the sum of the complexities of the blocks.
Recalling that in~\cite{Matveev:AAM} a complexity $c(M)$ was defined also
for $\partial M\ne\emptyset$, we emphasize here that our complexity
$c(M,X)$ is typically different from $c(M)$. So the splitting theorem crucially
depends on the extension of $c$ from manifolds to pairs.

Our splitting differs from the JSJ
decomposition (see \cite{Jaco-Shalen,Johannson}, and
\emph{e.g.}~\cite[Chapter 1]{Kapovich:book} for more recent developments)
for not being unique (see below for further discussion on this point),
but it has the great advantage that the blocks it involves, which we call \emph{bricks},
are much easier than all Seifert and
simple manifolds. As a matter of fact, our splitting is non-trivial on almost all
Seifert and hyperbolic manifolds it has been tested on.
Another advantage is that the graphs in the boundary reduce the
flexibility of possible gluings of bricks. As a consequence, a given set of bricks
can only be combined in a finite number of ways. This property is of course crucial
for computation, and our theory actually leads to very effective algorithms
for the enumeration of closed manifolds having small complexity.

Back to the relation of our splitting with the JSJ decomposition,
we mention that all the bricks found so far~\cite{MaPe:EM} are geometrically atoroidal,
which suggests that our splitting is actually always a refinement of the JSJ
decomposition. Moreover, non-uniqueness for a Seifert manifold
typically corresponds to non-uniqueness of its realization as a graph-manifold.
We also know of one non-uniqueness instance in the hyperbolic case.

The orientable version of the theory developed in this paper, culminating
in the splitting theorem, was established in~\cite{MaPe:EM}. In the same paper
we have proved several strong restrictions on the topology of bricks and, using
a computer program, we have been able to classify all orientable bricks of
complexity up to 9. Using the bricks we have then listed all closed irreducible
orientable 3-manifolds up to complexity 9, showing in particular that the
only four hyperbolic ones are precisely those of least known volume.
The splitting theorem proved below is the main theoretical tool needed to extend
our program of enumerating 3-manifolds of small complexity from the orientable to the
general case. We are planning to realize this program in the close future.
This will allow us to provide information on the smallest
non-orientable hyperbolic manifolds and
on the density, in each given complexity, of orientable
manifolds among all 3-manifolds.

We have decided to devote the present paper to the general theory
and the splitting theorem, leaving computer implementation for a subsequent paper,
because the non-orientable case displays certain remarkable phenomena
which do not appear in the orientable case. To begin with, toric boundary
components force the shape of the trivalent graph they contain to only
one possibility, while Klein bottles allow two. Next, the
assumption of $\matP^2$-irreducibility has to be added to
irreducibility to get finiteness of closed manifolds of a given complexity.
More surprisingly, these assumptions do not suffice when non-empty boundary
is allowed, because the drilling of a boundary-parallel orientation-reversing
loop never changes complexity. Because of these facts,
the intrinsic definition of \emph{brick} given below is somewhat subtler
than in~\cite{MaPe:EM}, and the proof of some of the key results (including
additivity under connected sum) is considerably harder.

\section{Manifolds with marked boundary}\label{manifolds:section}

If $C$ is a connected surface, we call \emph{spine} of $C$ a trivalent graph
$\tau$ embedded in $C$ in such a way that $C\setminus\tau$ is an open disc.
(A `graph' for us is just a `one-dimensional complex,' \emph{i.e.}~multiple
and closed edges are allowed.)
If $C$ is disconnected then a spine of $C$ is a collection of spines
for all its components.

We denote by $\calX$ the set of all pairs $(M,X)$, where $M$ is a connected and
compact 3-manifold with (possibly empty) boundary made of tori and Klein bottles,
and $X$ is a spine of $\partial M$.
Elements of $\calX$ will be viewed up to the natural equivalence relation
generated by homeomorphisms of manifolds.

\begin{rem}\label{all:zero:rem} \emph{If a connected surface $C$ has a
spine $\tau$ with $k\ge1$ vertices then $k$ is even and
$\chi(C)=k-3k/2+1\le0$. So, instead of specifying that for $(M,X)\in\calX$
the boundary $\partial M$ should consist of tori and Klein bottles, we may
have asked only that $\chi(M)$ should vanish and all elements of $X$
should have vertices.} \end{rem}

\begin{figure}\begin{center}
\figfatta{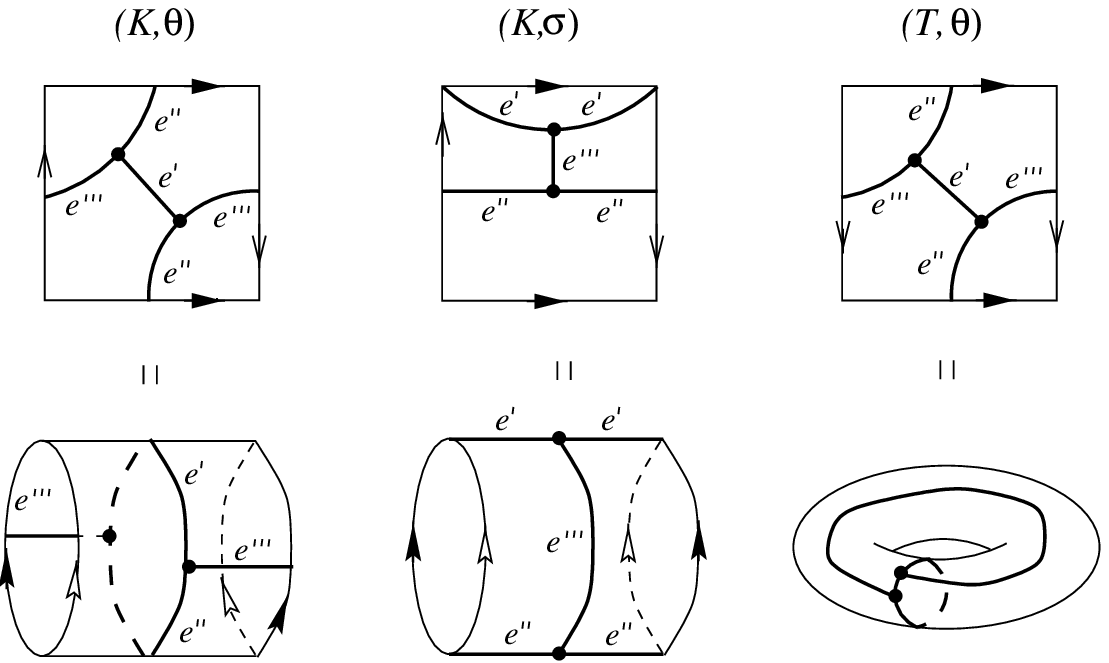}
\nota{Spines of the Klein bottle and the torus.}\label{twodimspin:fig}
\end{center}\end{figure}

\paragraph{Spines of the torus $T$ and the Klein bottle $K$}
A spine of $T$ or $K$ must be a trivalent graph with
two vertices, and there are precisely two such graphs, namely the $\theta$-curve
and the frame $\sigma$ of a pair of spectacles. Both $\theta$ and $\sigma$ can serve as spines
of the Klein bottle $K$, as suggested in Fig.~\ref{twodimspin:fig},
left and center.
The next result will be shown in the appendix:

\begin{prop}\label{Klein:spines:prop}
The following holds for both $\tau=\theta$ and $\tau=\sigma$:
\begin{enumerate}
\item The embedding of $\tau$ in $K$ described in
Fig.~\ref{twodimspin:fig} is the only one (up to isotopy)
such that $K\setminus\tau$ is an open disc.
\item There exists $f\in{\rm Aut}(K)$ such that $f(\tau)=\tau$ and $f$ interchanges
the edges $e'$ and $e''$, but every $f\in{\rm Aut}(K)$ such that $f(\tau)=\tau$
leaves $e'''$ invariant.
\end{enumerate}
\end{prop}

The situation for the torus $T$ is completely different. First of all, $\sigma$ is not
a spine of $T$. In addition, $\theta$ can be used as a spine of $T$ in infinitely many non-isotopic ways, because the position of $\theta$ on $T$ is determined by the triple of slopes
on $T$ which are contained in $\theta$. Note that these three slopes intersect each other
in a single point, and any such triple determines one spine $\theta$.
However we have the following result, which we leave to the reader to prove using
the facts just stated.

\begin{prop}\label{torus:spines:prop}
If $\theta$ is any spine of $T$ then all the
automorphisms of $\theta$ are induced by automorphisms of $T$.
If $\theta$ and $\theta'$ are spines of $T$ then there exists $f\in{\rm Aut}(T)$
such that $f(\theta)=\theta'$.
\end{prop}

\begin{figure}\begin{center}
\figfatta{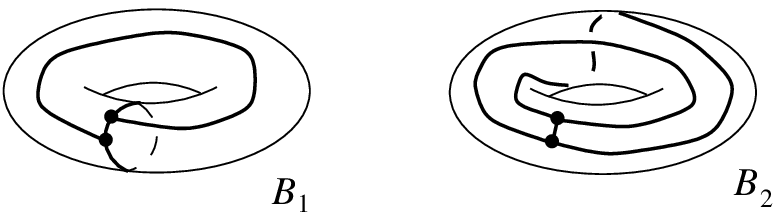}
\nota{The pairs $B_1$ and $B_2$.}\label{firstbricks:fig}
\end{center}\end{figure}

\paragraph{Examples of pairs}
Of course if $M$ is a closed 3-manifold then $(M,\emptyset)$ is an element of $\calX$.
For the sake of simplicity we will often write only $M$ instead of $(M,\emptyset)$.
We list here several more elements of $\calX$ which will be needed below.
Our notation will be consistent with that of~\cite{MaPe:EM}. The reader is invited to
use Propositions~\ref{Klein:spines:prop} and~\ref{torus:spines:prop} to make sure
that all the pairs we introduce are indeed well-defined up to homeomorphism.
We start with the product pairs:
\begin{eqnarray*}
B_0&=&(T\times[0,1],\{\theta\times\{0\},\theta\times\{1\}\}),\\
B'_0&=&(K\times[0,1],\{\theta\times\{0\},\theta\times\{1\}\}),\\
B''_0&=&(K\times[0,1],\{\sigma\times\{0\},\sigma\times\{1\}\}).
\end{eqnarray*}
We next have two pairs $B_1$ and $B_2$ based on the solid torus $\solT$
and shown in Fig.~\ref{firstbricks:fig}, and two on the solid
Klein bottle $\solK$, namely
$B'_1=(\solK,\theta)$ and $B'_2=(\solK,\sigma)$.

\begin{figure}\begin{center}
\figfatta{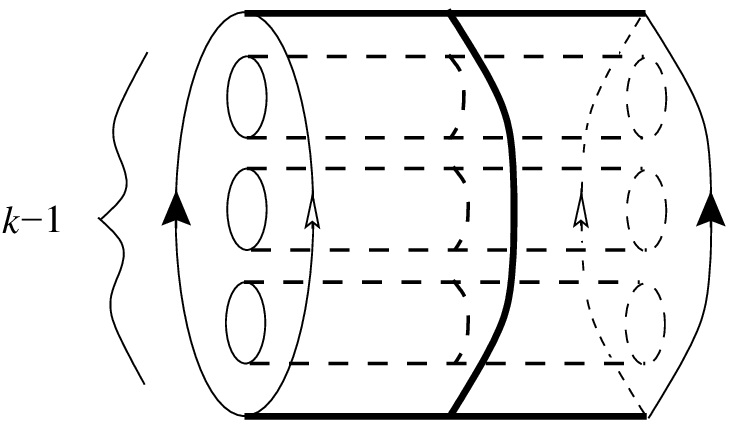}
\nota{The pair $Z_k$ for $k\ge1$.}\label{btwok:fig}
\end{center}\end{figure}

For $k\ge 1$ we consider now the 2-orbifold given by the disc $D^2$ with $k$ mirror
segments on $\partial D^2$. Then we define $Z_k\in\calX$ as the Seifert fibered
space without singular fibers over this 2-orbifold (see~\cite{Scott}),
with one spine $\sigma$ in each of the $k$ Klein bottles on the boundary.
Note that $Z_k$ can also be viewed as the complement of $k$
disjoint orientation-reversing loops in
$S^2\timtil S^1$. Yet another description of $Z_k$ is given in
Fig.~\ref{btwok:fig}.
We also note that $Z_1=B'_2$ and $Z_2=B''_0$. We define now $B''_2$ to be $Z_3$.
This notation has a specific reason explained below.

We will now introduce three operations on pairs which allow to construct
new pairs from given ones. The ultimate goal is to show that all manifolds
can be constructed via these operations using only certain
building blocks.

\paragraph{Connected sum of pairs} The operation of connected sum ``far from the
boundary'' obviously extends from manifolds to pairs.
Namely, given $(M,X)$ and $(M',X')$ in $\calX$, we define $(M,X)\#(M',X')$
as $(M\# M',X\cup X')$, where $M\# M'$ is one of the two possible
connected sums of $M$ and $M'$. Of course $S^3=(S^3,\emptyset)\in\calX$ is the identity
element for operation $\#$. It is now natural to define $(M,X)$ to be \emph{prime} or
\emph{irreducible} if $M$ is. Of course the only prime non-irreducible pairs
are $S^2\times S^1=(S^2\times S^1,\emptyset)$ and $S^2\timtil S^1=(S^2\timtil S^1,\emptyset)$.

\paragraph{Assembling of pairs} Given $(M,X)$ and $(M',X')$ in $\calX$, we pick
spines $\tau \in X$ and $\tau' \in X'$ of the same type $\theta$ or $\sigma$.
If $\tau\subset C\subset\partial M$ and $\tau'\subset C'\subset\partial M'$
we choose now a
homeomorphism $\psi:C \to C'$ such that $\psi(\tau) =
\tau'$. We can then construct the pair
$(N,Y) = (M \cup_{\psi} M',(X\cup X')\setminus \{\tau, \tau'\})$.
We call
this an \emph{assembling} of $(M,X)$ and $(M',X')$ and we write
$(N,Y)=(M,X)\oplus(M',X')$. Of course two given elements of $\calX$ can
only be assembled in a finite number of inequivalent ways.

Considering the pairs $B_i^*$ and $Z_k$ introduced above,
the reader may easily check as an exercise that $Z_k\oplus Z_h=Z_{h+k-2}$
and that the following holds:

\begin{rem}\label{trivial:assemblings:rem}
\begin{enumerate}
\item[{\rm 1.}] \emph{$(M,X)\oplus B_0^*=(M,X)$ for any $(M,X)\in\calX$;}
\item[{\rm 2.}] \emph{It is possible to assemble the pair $B_1$ to itself along a certain map $\psi$
so to get $\closedpair{S^3}$. This implies that, starting from $(M,X)$, if we first perform a
connected sum $(M,X)\# B_1$ and then the assembling $((M,X)\# B_1)\oplus B_1$ along this map
$\psi$, we get the original $(M,X)$ as a result. Similarly one can assemble $B_2$ and $B_1$
so to get $\closedpair{S^3}$, whence $((M,X)\# B_2)\oplus B_1=
((M,X)\# B_1)\oplus B_2=(M,X)$ (for suitably chosen gluings);}
\item[{\rm 3.}] \emph{The assembling of $B''_2$ with $B'_2$ gives
$B''_0$, so $((M,X)\oplus B''_2)\oplus B'_2=(M,X)$
provided $B'_2$ is assembled to one of the free boundary components of $B''_2$.}
\end{enumerate}
\end{rem}

\noindent This remark shows that we can discard various assemblings without impairing
our capacity of constructing new manifolds. To be precise we will call
\emph{trivial} an assembling $(M,X)\oplus (M',X')$ if, up to
interchanging $(M,X)$ and $(M',X')$, one of the following holds:
\begin{enumerate}
\item $(M',X')$ is of type $B^*_0$;
\item $(M',X')=B_j$ for $j\in\{1,2\}$ and
$(M,X)$ can be expressed as $(N,Y)\# B_i$ for $i\in\{1,2\}$ with $(N,Y)\ne\closedpair{S^3}$
in such a way that the assembling is performed along the boundary of $B_i$ and
$B_i\oplus B_j=\closedpair{S^3}$;
\item $(M',X')=B'_2$ and $(M,X)$ can be expressed as $(N,Y)\oplus B''_2$
with $B'_2$ being assembled to $B''_2$.
\end{enumerate}

\paragraph{Self-assembling} Given $(M,X)\in\calX$, we pick two distinct
spines $\tau, \tau'\in X$ with $\tau\subset C$
and $\tau'\subset C'$. We choose a homeomorphism $\psi:C
\to C'$ such that $\psi(\tau)$ and $\tau'$ intersect
transversely in two points, and we construct the pair
$(N,Y) = (M_{\psi}, X \setminus \{\tau, \tau'\})$. We call this
a \emph{self-assembling} of $(M,X)$ and we write $(N,Y) = \odot(M,X)$. As
above, only a finite number of self-assemblings of a given element of
$\calX$ are possible.

In the sequel it will be convenient to refer to a combination of
assemblings and self-assemblings of pairs just as an \emph{assembling}.
Note that of course we can do the assemblings first and the
self-assemblings in the end.

\section{Complexity, bricks, and the decomposition theorem}\label{complexity:section}
Starting from the next section we will introduce and discuss
a certain function $c:\calX\to\matN$ which we call \emph{complexity}.
In the present section we only very briefly anticipate the definition of $c$
and state several results about it, which could also be taken as axiomatic properties.
Then we show how to deduce the splitting theorem from the properties only.
Proofs of the properties are given in Sections~\ref{skeleta:section}
to~\ref{assemblings:section}.

Given $(M,X)\in\calX$ we denote by $c(M,X)$ and call the \emph{complexity} of $(M,X)$ the
minimal number of vertices of a simple polyhedron $P$ embedded in $M$ such that
$P\cup\partial M$ is also simple, $P\cap\partial M=X$, and
the complement of $P\cup\partial M$ is an open 3-ball.
Here `simple' means that the link of every point embeds in the 1-skeleton of
the tetrahedron, and a point of $P$ is a `vertex' if its link is precisely the 1-skeleton of
the tetrahedron. We obviously have:

\begin{prop}\label{old:complexity:prop}
If $M$ is a closed 3-manifold then $c(M)=c(M,\emptyset)$
coincides with
Matveev's $c(M)$ defined in~\cite{Matveev:AAM}.
\end{prop}

Note that $c(M)$ is also defined in~\cite{Matveev:AAM}
for $\partial M\ne\emptyset$, but
typically $c(M,X)\ne c(M)$.

\paragraph{Axiomatic properties}
We start with three theorems which suggest to restrict the study of $c(M,X)$
to pairs $(M,X)$ which are irreducible and \ptwoirred. Recall that $M$
is called \emph{\ptwoirred}\ if it does not contain any two-sided embedded
projective plane $\matP^2$ (see~\cite{Hempel:book} for generalities about
this notion, in particular for the proof that
a connected sum is \ptwoirred\ if and only if the individual summands are).
When $M$ is closed, we call \emph{singular} a triangulation of $M$ with multiple
and self-adjacencies between tetrahedra.
The first and second theorems extend results of Matveev~\cite{Matveev:AAM}
respectively from the closed to the marked-boundary case, and
from the orientable to the possibly-non-orientable case.
The extension is easy for the second theorem, not quite so
for the first theorem. The third theorem shows that the
non-orientable theory is far richer than the orientable one.

\begin{teo}[additivity under $\#$]\label{additivity:teo}
For any $(M,X)$ and $(M',X')$ we have
$$c((M,X)\#(M',X'))=c(M,X)+c(M,X').$$
Moreover $c(S^2\times S^1)=c(S^2\timtil S^1)=0$.
\end{teo}

\begin{teo}[naturality]\label{naturality:teo}
If $M$ is closed, irreducible, \ptwoirred, and different from
$S^3$, $\matP^3$, $L_{3,1}$,
then $c(M)=c(M,\emptyset)$ is the minimal number of tetrahedra in a singular
triangulation of $M$.
\end{teo}

\begin{teo}[finiteness]\label{finiteness:teo}
For all $n\ge0$ the following happens:
\begin{enumerate}
\item There exist finitely many irreducible and \ptwoirred\ pairs $(M,X)$ such
that $c(M,X)=n$ and
$(M,X)$ cannot be expressed as an assembling $(N,Y)\oplus B''_2$;
\item If $(N,Y)\in\calX$ is irreducible and \ptwoirred\
and $c(N,Y)=n$ then $(N,Y)$ can be obtained from
one of the $(M,X)$ described above by repeated assembling of copies of $B''_2$.
Any such assembling has complexity $n$.
\end{enumerate}
\end{teo}

The previous result is of course crucial for computational purposes.
To better appreciate its ``finiteness'' content, note that
whenever we assemble one copy of $B''_2$ the number of boundary components
increases by one. Therefore the theorem implies that for all $n,k\ge 0$ the set
$$\calM_{\le n}^{\le k}=\{(M,X)\in\calX\
{\rm irred.\ and}\ \matP^2{\rm -irred.},\ c(M,X)\le n,\ \# X\le k\}$$
is finite. It should be emphasized that not only can we
prove that $\calM_{\le n}^{\le k}$ is finite, but the proof itself
provides an explicit algorithm to produce a finite list of pairs from which
$\calM_{\le n}^{\le k}$ is obtained by removing duplicates. The theorem
also implies that dropping the restriction $\#X\le k$ we get infinitely many pairs,
but only finitely many orientable ones.
This fact, which is ultimately due to the existence of the $Z_k$ series
generated by $B''_2$ under assembling, is one of the key differences between
the orientable and the general case
(another important difference will arise in the proof of
Theorem~\ref{additivity:teo} ---see Proposition~\ref{normal:sphere:prop}).
Note also that an assembling with $B''_2$ geometrically corresponds
to the drilling of a boundary-parallel orientation-reversing loop.
A more specific version of the previous theorem for $n=0$ is needed below:

\begin{prop}\label{compl:zero:prop}
The only irreducible and \ptwoirred\ pairs having complexity $0$ are
$\closedpair{S^3}$, $\closedpair{L_{3,1}}$, $\closedpair{\matP^3}$ and
all the $B_i^*$ and $Z_k$ defined above.
\end{prop}

We turn now to the behavior of complexity under
assembling. All the results stated in the rest of this section are
new and strictly depend on the extension to pairs of the theory
of complexity.

\begin{prop}[subadditivity]\label{subadditivity:prop}
For any $(M,X),(M',X')\in\calX$ we have:
\begin{eqnarray*}
c((M,X)\oplus(M',X'))&\le& c(M,X)+c(M',X'),\\
c(\odot(M,X))&\le& c(M,X)+6.
\end{eqnarray*}
\end{prop}

We define now an assembling $(M,X)\oplus(M',X')$ to be \emph{sharp}
if it is non-trivial and $c((M,X)\oplus(M',X'))=c(M,X)+c(M',X')$.
Similarly, a self-assembling $\odot(M,X)$ is \emph{sharp} if
$c(\odot(M,X))= c(M,X)+6$.
Proposition~\ref{subadditivity:prop} readily implies the following:

\begin{rem}\label{B2sec:is:sharp:rem}
\begin{enumerate}
\item[{\rm 1.}] \emph{If a combination of sharp
(self-)assemblings is rearranged in a different order then
it still consists of sharp (self-)assemblings;}
\item[{\rm 2.}] \emph{Every assembling with $B''_2$ is sharp (unless it is trivial, which only happens
when $B''_2$ is assembled to $B''_0$ or to $B'_2$).
To see this, note again that $(M,X)\oplus B''_2\oplus B'_2=(M,X)$ and
$c(B''_2)=c(B'_2)=0$.}
\end{enumerate}
\end{rem}

\begin{teo}[sharp splitting]\label{sharp:split:teo}
Let $(N,Y)$ be irreducible and \ptwoirred. If $(N,Y)$ can be expressed as a sharp
assembling $(M,X)\oplus(M',X')$ or as a self-assembling
$\odot(M'',X'')$ then $(M,X)$, $(M',X')$, and $(M'',X'')$ are irreducible and
\ptwoirred.
\end{teo}

\begin{proof}
In both cases we are cutting $N$ along a two-sided torus or Klein bottle,
so $\matP^2$-irreducibility is obvious. If $(N,Y)=\odot(M'',X'')$, this torus
or Klein bottle is incompressible in $N$, and irreducibility of $M''$ is a general
fact~\cite{Hempel:book}. We are left to show that if $(N,Y)=(M,X)\oplus(M',X')$ sharply
then $M$ and $M'$ are irreducible. Since they have boundary, it is enough to
show that they are prime. Suppose they are not, and
consider prime decompositions of $(M,X)$ and $(M',X')$
involving summands $(M_i,X_i)$ and $(M'_j,X'_j)$. So one summand
$(M_i,X_i)$ is assembled to one $(M'_j,X'_j)$, and the other
$(M_i,X_i)$'s and $(M'_j,X'_j)$'s survive in $(N,Y)$. It follows that, up to permutation,
$(M,X)$ is prime, $(M',X')=(M'_1,X'_1)\#(M'_2,X'_2)$ with
$(M'_1,X'_1)$ and $(M'_2,X'_2)$ prime, $(M,X)\oplus (M'_1,X'_1)=
\closedpair{S^3}$ and $(M'_2,X'_2) = (N,Y)$.
Sharpness of the original assembling and additivity under $\#$ now imply that $c(M,X)=c(M'_1,X'_1)=0$.
So Proposition~\ref{compl:zero:prop} applies to $(M,X)$ and $(M'_1,X'_1)$.
Knowing that $(M,X)\oplus (M',X')=\closedpair{S^3}$
it is easy to deduce that $(M,X)$ and $(M'_1,X'_1)$ are either $B_1$ or $B_2$,
and that the original assembling was a trivial one. A contradiction.
\end{proof}

\paragraph{Bricks and decomposition}
Taking the results stated above for granted, we define here the elementary
building blocks and prove the decomposition theorem. Later we will make comments
about the actual relevance of this theorem.

A pair $(M,X)\in\calX$ is called a \emph{brick} if it is irreducible and
\ptwoirred\ and cannot be expressed as a sharp assembling or self-assembling.
Theorem~\ref{finiteness:teo} and Remark~\ref{B2sec:is:sharp:rem} easily
imply that there are finitely many bricks of complexity $n$.
From Proposition~\ref{compl:zero:prop} it is easy to deduce that in complexity
zero the only bricks are precisely the $B_i^*$ introduced above,
which explains why we have given a special status to $Z_3=B''_2$,
and that the other irreducible and \ptwoirred\ pairs are assemblings of bricks.
Now, more generally:

\begin{teo}[existence of splitting]\label{splitting:teo}
Every irreducible and \ptwoirred\ pair $(M,X)\in\calX$ can be expressed as a
sharp assembling of bricks.
\end{teo}

\begin{proof}
The result is true for $c(M,X)=0$, so
we proceed by induction on $c(M,X)$ and suppose $c(M,X)>0$.
By Theorem~\ref{finiteness:teo} we can assume that $(M,X)$
cannot be split as $(N,Y)\oplus B''_2$, because
every assembling with $B''_2$ is sharp, and we
have seen that $B''_2$ is a brick. Now if $(M,X)$ is a brick we are done.
Otherwise $(M,X)$ is either a sharp self-assembling $\odot(N,Y)$, but in this case
$c(N,Y)=c(M,X)-6$ and we conclude by induction using Theorem~\ref{sharp:split:teo}, or
$(M,X)$ is a sharp assembling $(N,Y)\oplus(N',Y')$. Theorem~\ref{sharp:split:teo}
states that $(N,Y)$ and $(N',Y')$ are irreducible and \ptwoirred.
If both $(N,Y)$ and $(N',Y')$
have positive complexity we conclude by induction. Otherwise we can assume that
$c(N',Y')=0$ and apply Proposition~\ref{compl:zero:prop}.
Since the assembling is non-trivial, $(N',Y')$ is not
of type $B_0^*$. It is also not $B''_2$ or $Z_k$ for $k\ge 3$,
by the property of $(M,X)$ we
are assuming. So $(N',Y')$ is one of $B_1$, $B'_1$, $B_2$, $B'_2$.
In particular, it is a brick.

Now we claim that $(N,Y)$ cannot be split as $(N'',Y'')\oplus B''_2$.
Assuming it can, we have two cases.
In the first case the assembling of $(N',Y')$ is performed along a free
boundary component of $B''_2$, but then we must have $(N',Y')=B'_2$, and
the assembling is trivial, which is absurd. In the second case
$(N',Y')$ is assembled to a free boundary component of $(N'',Y'')$, and we have
$$(M,X)=\big((N'',Y'')\oplus(N',Y')\big)\oplus B''_2,$$
which is again absurd. Our claim is proved.

Now we know that $(N,Y)$ again belongs to the finite list of
irreducible and \ptwoirred\ manifolds which have complexity $n$
and cannot be split as an assembling with $B''_2$. However $(N,Y)$
has one more boundary component than $(M,X)$, which implies that by
repeatedly applying this argument we must eventually end up with a brick.
\end{proof}

\paragraph{Classification of bricks}
Theorem~\ref{splitting:teo} shows that listing irreducible
and \ptwoirred\ manifold up to complexity
$n$ is easy once the bricks up to complexity
$n$ are classified. The finiteness features of our theory imply
that there exists an algorithm which reduces such a classification
to a recognition problem. We illustrate here this algorithm
and give a hint to explain why does it work in practice. To do this we will
need to refer to results stated and proved later in the paper.

We know the bricks of complexity zero, so we fix $n\ge 1$ and
inductively assume to know the set $\calB_{<n}$ of bricks
of complexity up to $n-1$. Theorem~\ref{standard:teo}
implies that there exists an effective method
to produce a finite list $\calL_n$ which contains (with repetitions)
all irreducible and \ptwoirred\ pairs $(M,X)$
such that $c(M,X)\le n$ and $\#X\le 2n$, and
Corollary~\ref{limitazione:su:X:cor} now implies that all bricks of
complexity $n$ appear in the list.

Suppose now that for some reason we can extract from $\calL_n$ a shorter
list $\calL'_n$ which we know to still contain all bricks
of complexity $n$. We also assume that $\calL'_n$ does not contain pairs of
complexity zero. To make sure that a given element $(M,X)$ of $\calL'_n$
is a brick we must now check that it is not homeomorphic to a
sharp assembling of elements of $\calB_{<n}$ and other elements of
$\calL'_n$. In a sharp assembling of bricks giving $(M,X)$ we can of
course have at most $n$ positive-complexity bricks,
and the knowledge of the bricks of complexity zero
shows that we can also have at most $2n$ bricks of complexity zero.
Therefore, to check whether $(M,X)$ is a brick, we only need to recognize
whether it belongs to a finite list of pairs.

Besides the recognition problem, the crucial step of the algorithm
just described is the extraction of the list $\calL'_n$ from
the list $\calL_n$. The point is that $\calL_n$ is hopelessly big even
for small $n$, so to actually classify bricks one must be able to produce
a much shorter $\calL'_n$ without even knowing the whole of $\calL_n$.
This was achieved in~\cite{MaPe:EM}, in the orientable case with
$n\le 9$, by means of a number of results which provide
strong \emph{a priori} restrictions on the topology of the bricks.
As explained in the introduction, the
non-orientable version of these results and the computer search of
the first non-orientable bricks are deferred to a subsequent paper.

\paragraph{Interesting assemblings}
The practical relevance of Theorem~\ref{splitting:teo} towards the classification
of irreducible and \ptwoirred\ 3-manifolds of bounded complexity sits in the following
heuristic facts:
\begin{enumerate}
\item For any $n$ the number of bricks of complexity at most $n$ is by far
smaller than the number of all irreducible and \ptwoirred\ pairs, and
the above-described algorithm to find the bricks is rather efficient;
\item If a manifold is expressed as an assembling of bricks, it is typically
easy to recognize the manifold and its JSJ decomposition, and hence to make
sure that the assembling is sharp by checking that the same manifold was not obtained
already in lower complexity;
\item When an assembling of bricks
is sharp, it is typically true that the result
is again irreducible and \ptwoirred.
\end{enumerate}
Facts 1 and 2 can be made precise when $n\le 9$ and only orientable manifolds are considered.
Namely it was shown in~\cite{MaPe:EM} that:
\begin{enumerate}
\item There are 1902 closed, irreducible, and orientable 3-manifolds of complexity up to 9,
and only 7 bricks can be used to obtain all but 19 of them.
(The other 19 manifolds are themselves bricks, but since they have empty
boundary they cannot be assembled at all.)
\item All the orientable bricks up to complexity 9 are geometrically atoroidal, so,
for a closed orientable $M$ with $c(M)\le 9$, each block of
the JSJ decomposition of $M$ is a union of some of the bricks of our decomposition.
\end{enumerate}

Concerning fact 3, we make it more precise here for both the orientable and
the non-orientable case.

\begin{teo}\label{sharp:assemble:teo}
\begin{enumerate}
\item Assume
$(M,X)$ and $(M',X')$ are irreducible and $\matP^2$-irredu\-cible pairs and
$(N,Y)=(M,X)\oplus (M',X')$ is a sharp assembling. Then $(N,Y)$ is
prime. It can fail to be \ptwoirred\ only if one of $M$ or $M'$ is a solid torus or a solid
Klein bottle.
\item Assume $(M'',X'')$ is irreducible and \ptwoirred\
and $(N,Y)=\odot(M'',X'')$ is a self-assembling. Then $(N,Y)$ is irreducible
and \ptwoirred.
\end{enumerate}
\end{teo}

\section{Skeleta}\label{skeleta:section}

In this section we introduce the notion of \emph{skeleton} of a pair $(M,X)$,
we define the complexity of $(M,X)$ as the minimal number of vertices of a skeleton,
and we discuss the first properties of minimal skeleta, deducing some of the
results stated above. The other results, which require a deeper analysis and
new techniques, will be proved in subsequent sections.

\paragraph{Simple skeleta and definition of complexity}
We recall that a compact polyhedron $P$ is called \emph{simple} if the link of every
point of $P$ can be embedded in the space given by a circle with three
radii. The points having the whole of this space as a link are called
\emph{vertices}. They are isolated and therefore finite in number.

Given a pair $(M,X)\in\calX$, a polyhedron $P$ embedded in $M$ is called
a \emph{skeleton} of $(M,X)$ if the following conditions hold:
\begin{itemize}
\item $P\cup \partial M$ is simple;
\item $M\setminus (P\cup\partial M)$ is an open ball;
\item $P\cap\partial M = X$.
\end{itemize}

\begin{rem}\label{skeleton:is:simple:rem}
\emph{If $P$ is a skeleton of $(M,X)$ then $P$ is simple, and
the vertices of $P$ cannot lie on $\partial M$.
When $\#X = 1$ then $P$ is a \emph{spine} of $M$ (\emph{i.e.}~$M$ collapses onto
$P$), and when $\#X=0$ (\emph{i.e.}~when $M$ is closed) then $P$ is a
spine in the usual sense~\cite{Matveev:AAM}, namely $M\setminus\{ {\rm point } \}$
collapses onto $P$. When $\#X \ge 2$ no such interpretation is possible.}
\end{rem}

\begin{rem} \emph{It is easy to prove that every $(M,X)\in\calX$ has a
skeleton: take any simple spine $Q$ of $M\setminus\{{\rm point}\}$, so that
$M\setminus Q=\partial M \times [0,1)\cup B^3$, and assume that,
as $\tau$ varies in $X$, the
various $\tau\times[0,1)$'s are incident in a generic way to $Q$
and to each other. Taking the union of $Q$ with the $\tau \times [0,1)$'s we
get a simple $Q'$ such that $M\setminus(Q'\cup\partial M)$ consists
of $\#X+1$ balls. Then we get a skeleton of $(M,X)$ by puncturing
$\#X$ suitably chosen 2-discs embedded in $Q'$, so to get one ball
only in the complement.}
\end{rem}

\begin{rem}{\em A definition of skeleton analogous to our one was given
in~\cite{Tu-Vi} for any compact manifold with any trivalent graph in its
boundary.} \end{rem}

For a simple polyhedron $P$ we denote by $v(P)$ the number of vertices of $P$, and
we define the \emph{complexity} $c(M,X)$ of a given $(M,X)\in\calX$ as the minimum of
$v(P)$ over all skeleta $P$ of $(M,X)$. So we have a function
$c:\calX\to\matN$.

\paragraph{Some skeleta without vertices}
If we remove one point from the closed manifolds $S^3$, $L_{3,1}$, $\matP^3$,
$S^2\times S^1$, and $S^2\timtil S^1$ then we can collapse the result
respectively to a point, to the ``triple hat,'' to the projective plane, and to the join
of $S^2$ and $S^1$ (for both the last two cases). Here the triple hat is the space obtained
by attaching the disc to the circle so that the boundary of the disc runs
three times around the circle. This shows that
$S^3$, $L_{3,1}$, $\matP^3$,
$S^2\times S^1$, and $S^2\timtil S^1$ all have complexity zero. It is a well-known fact,
which we will prove again below, that these are the only prime and \ptwoirred\
manifolds having complexity zero.

Turning to the $B_i^*$ and $Z_k$ defined in the previous section, we now show that they also have complexity $0$. This is rather obvious for the product pairs $B_0$, $B'_0$, and
$B''_0$, because they have the product skeleta $P_0=\theta\times[0,1]\subset T\times[0,1]$,
$P'_0=\theta\times[0,1]\subset K\times[0,1]$,
and $P''_0=\sigma\times[0,1]\subset K\times[0,1]$.

\begin{figure}\begin{center}
\figfatta{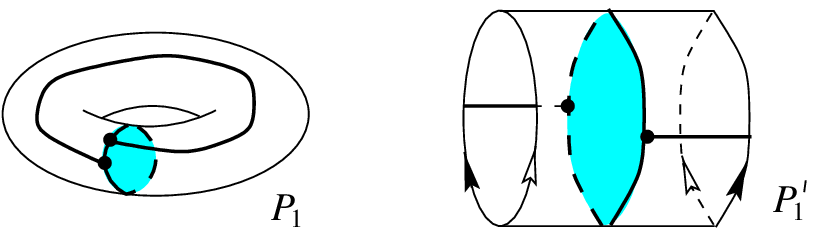}
\nota{The skeleta $P_1$ and $P'_1$ of $B_1$ and $B'_1$.}\label{p1spines:fig}
\end{center}\end{figure}

For $B_1=(\solT,\{\theta\})$
we note that $\theta$ contains a meridian of the torus, so we can attach
to $X$ a meridional disc and get the skeleton $P_1$ shown in Fig.~\ref{p1spines:fig}.
The same construction applies to $B'_1=(\solK,\{\theta\})$ and leads to the skeleton
$P'_1$ also shown in the figure.
Of course $P_1$ and $P'_1$ are isomorphic as abstract polyhedra
(just as $P_0$ and $P'_0$), but
we use different names to keep track also of their embeddings.

\begin{figure}\begin{center}
\figfatta{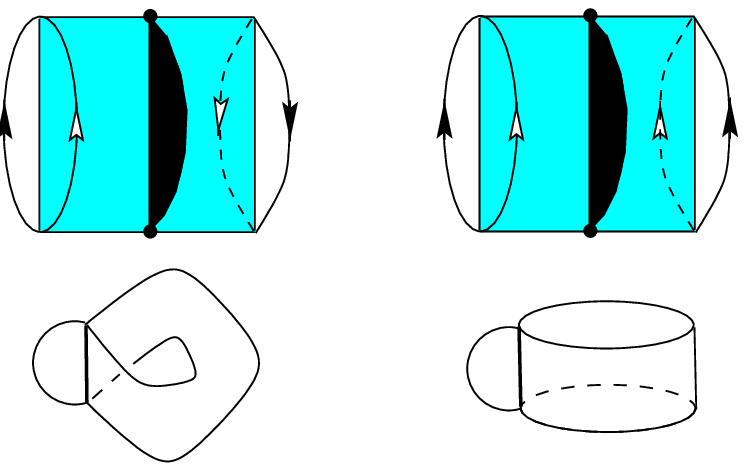}
\nota{The skeleta $P_2$ and $P'_2$ of $B_2$ and $B'_2$.}\label{p2spines:fig}
\end{center}\end{figure}

\begin{figure}\begin{center}
\figfatta{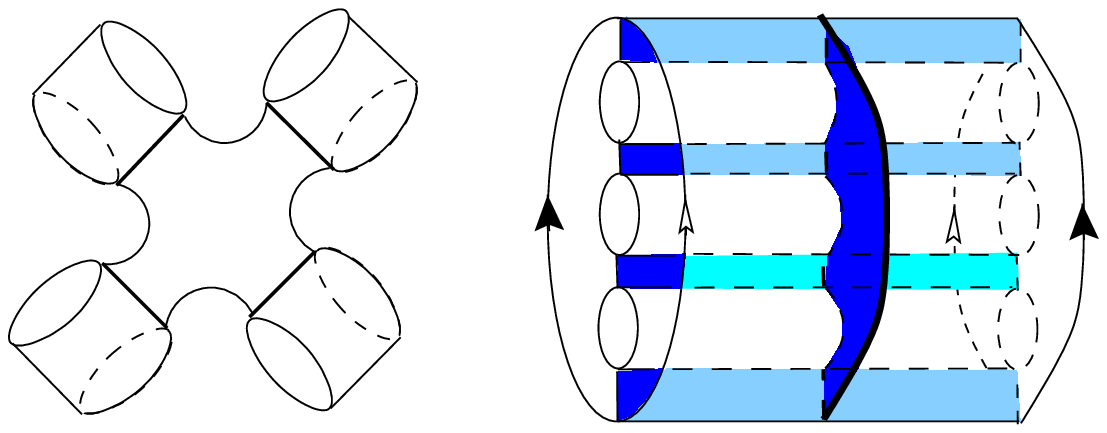}
\nota{The skeleton of $Z_k$ for $k=4$.}\label{p2kspines:fig}
\end{center}\end{figure}

Skeleta $P_2$ and $P'_2$ of $B_2$ and $B'_2$ respectively
are shown in Fig.~\ref{p2spines:fig}, both as abstract polyhedra and as embedded in $\solT$ and $\solK$.
We conclude with the series $Z_k$ for $k\ge 3$, for which a skeleton is
shown in Fig.~\ref{p2kspines:fig}. Recalling that $B''_2$ was defined as $Z_3$,
we denote this skeleton by $P''_2$ when $k=3$.

\paragraph{Nuclear, quasi-standard, and standard skeleta}
A skeleton of $(M,X)$ is called \emph{nuclear} if it does not
collapse to a proper subpolyhedron which is also a skeleton of $(M,X)$.
A nuclear skeleton $P$ of $(M,X)\in\calX$ having $c(M,X)$ vertices
is called \emph{minimal}. Of course every $(M,X)$ has minimal skeleta.

\begin{figure}
\begin{center}
\figfatta{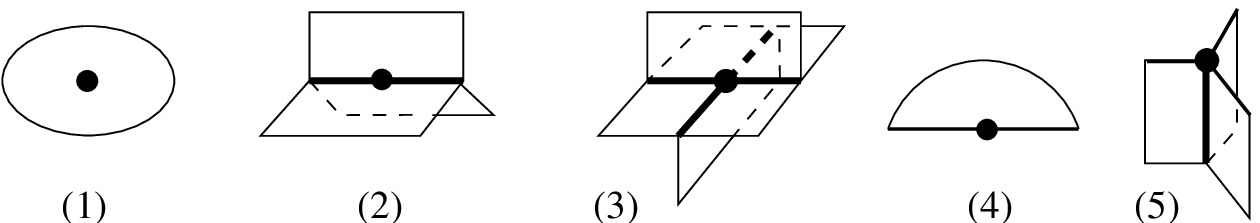}
\nota{Typical neighborhoods of points in a quasi-standard
polyhedron with boundary.} \label{quasistandard:fig}
\end{center}
\end{figure}

We will introduce now two more restricted classes of simple polyhedra.
Later we will show that, under suitable assumptions,
minimal polyhedra must belong to these classes.
A simple polyhedron $Q$ is called \emph{quasi-standard with boundary} if every point
has a neighborhood of one of the types (1)-(5) shown in Fig.~\ref{quasistandard:fig}.
A point of type (3) was already defined above to be a \emph{vertex} of $Q$. We denote
now by $V(Q)$ the set of all vertices, and we define the
\emph{singular set} $S(Q)$ as the set of points of type (2), (3), or (5), and
the \emph{boundary} $\partial Q$ as the set of points of type (4) or (5).
Moreover we call \emph{$1$-components} of $Q$ the connected components of $S(Q)
\setminus V(Q)$ and \emph{$2$-components} of $Q$ the connected components
of $Q \setminus (S(Q)\cup\partial Q)$.

If the $2$-components of $Q$ are open discs (and
hence are called just \emph{faces}), and the $1$-components are open
segments (and hence called just {\em edges}), then we call $Q$ a
\emph{standard polyhedron with boundary}. For short we will often just
call $Q$ a \emph{standard} polyhedron, and possibly specify that $\partial
Q$ should or not be empty.
We prove now the first properties of nuclear skeleta.

\begin{lemma} \label{nuclear:lem}
If $P$ is a nuclear skeleton of a pair $(M,X)\in\calX$,
then $P=Q\cup s_1\cup\ldots\cup s_m\cup G$, where:
\begin{enumerate}
\item $Q$ is a quasi-standard polyhedron with boundary $\partial Q \subset X$;
\item For all components $(C,\tau)$ of $(\partial M,X)$, either
$\partial Q\supset\tau$ or $Q$ appears near $C$ as in
Fig.~\ref{filialbordo:fig}, so $\partial Q\cap\tau$ is
one or two circles, depending on the type of $(C,\tau)$;
\begin{figure}\begin{center}
\figfatta{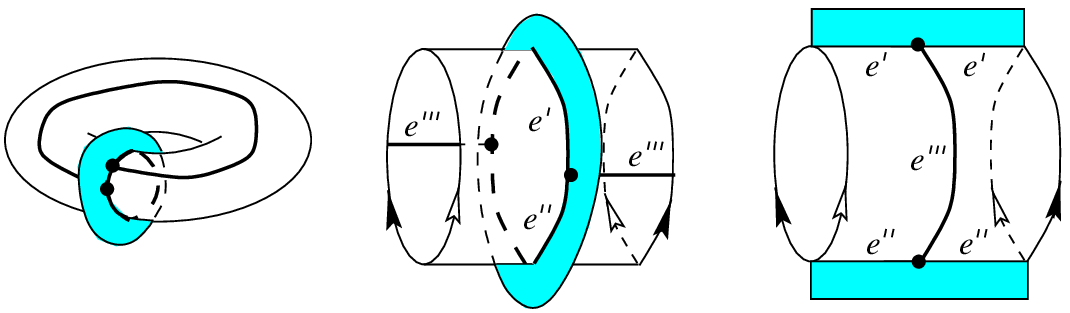}
\nota{Local aspect of $Q$ near $C$ if $\partial Q\not\supset\tau$.}\label{filialbordo:fig}
\end{center}\end{figure}
\item $s_1,\ldots,s_m$ are the edges of the $\tau$'s in $X$ which do not
already belong to $Q$;
\item $G$ is a graph with $G\cap (Q\cup s_1\cup\ldots\cup s_m)$ finite and
$G\cap V(Q\cup\partial M)$ empty.
\end{enumerate}
\end{lemma}

\begin{proof}
Nuclearity is a property of local nature, and the result is trivial if
$\partial M=\emptyset$. For $\partial M\ne\emptyset$, defining
$Q$ as the 2-dimensional portion of $P$
and $G$ as $P\setminus(Q\cup X)$, the only non-obvious point to show is (2).
Of course $Q\cap C\subset \partial Q$ is either $\tau$ or a union
of circles. To check that the only possibilities are those
of Fig.~\ref{filialbordo:fig} one recalls that
$M\setminus(P\cup\partial M)$ is a ball, so $C\setminus (Q\cup G)$
is planar, and then $C\setminus Q$ is also planar.
\end{proof}

\begin{rem}\label{nuclear:rem}
\emph{Every $(M,X)\in\calX$ has a
minimal skeleton $P'=Q\cup s_1\cup\ldots\cup s_m\cup G'$ as above,
where in addition $G'\cap\partial M=\emptyset$. This is because,
without changing $v(P)$, we can take the ends of $G$ lying on
$\partial M$ and make them slide over $Q\cup s_1\cup\ldots\cup s_m$
until they reach $\interior{M}$. Note that the regular neighborhood
of $\tau \in X$ in $P'$ is now either a product $\tau\times
[0,1]$ or as shown in Fig.~\ref{filialbordo:fig}.}
\end{rem}

\paragraph{Subadditivity}
Some properties of complexity readily follow from the definition and
from the first facts shown about minimal skeleta. To begin with, if $P$
and $P'$ are skeleta of $(M,X)$ and $(M',X')$ and we add to $P\sqcup P'$
a segment which joins $P\setminus V(P)$ to
$P'\setminus V(P')$, we get a skeleton of $(M,X)\#(M',X')$ with
$v(P)+v(P')$ vertices. This shows that $c((M,X)\#(M',X'))\le
c(M,X)+c(M',X')$.
Turning to assembling, let $P$
and $P'$ be minimal skeleta of $(M,X)$ and $(M',X')$ as in
Remark~\ref{nuclear:rem}, and let an assembling
$(M,X)\oplus(M',X')$ be performed along a map $\psi:C\to C'$ with
$\psi(\tau)=\tau'$. Then $P\cup_\psi P'$ is simple, and it is a skeleton of
$(M,X)\oplus(M',X')$. We deduce that
$c((M,X)\oplus(M',X'))\le c(M,X)+c(M',X')$.

Now we consider a self-assembling $\odot(M,X)$. If $P$ is a skeleton
of $(M,X)$ as in Remark~\ref{nuclear:rem}
and the self-assembling is performed along a certain map
$\psi:C\to C'$ such that $\tau'\cap\psi(\tau)$ consists of two points,
then $(P\cup C\cup C')/_\psi$ is a skeleton of $\odot(M,X)$. It has the same vertices
as $(M,X)$ plus at most two from the vertices of $\tau$,
two from the vertices of $\tau'$,
and two from $\tau'\cap\psi(\tau)$. This shows
that $c(\odot(M,X))\le c(M,X)+6$.

\paragraph{Surfaces determined by graphs}
We will need very soon the idea of splitting a skeleton along a graph,
so we spell out how the construction goes.
\begin{figure}\begin{center}
\figfatta{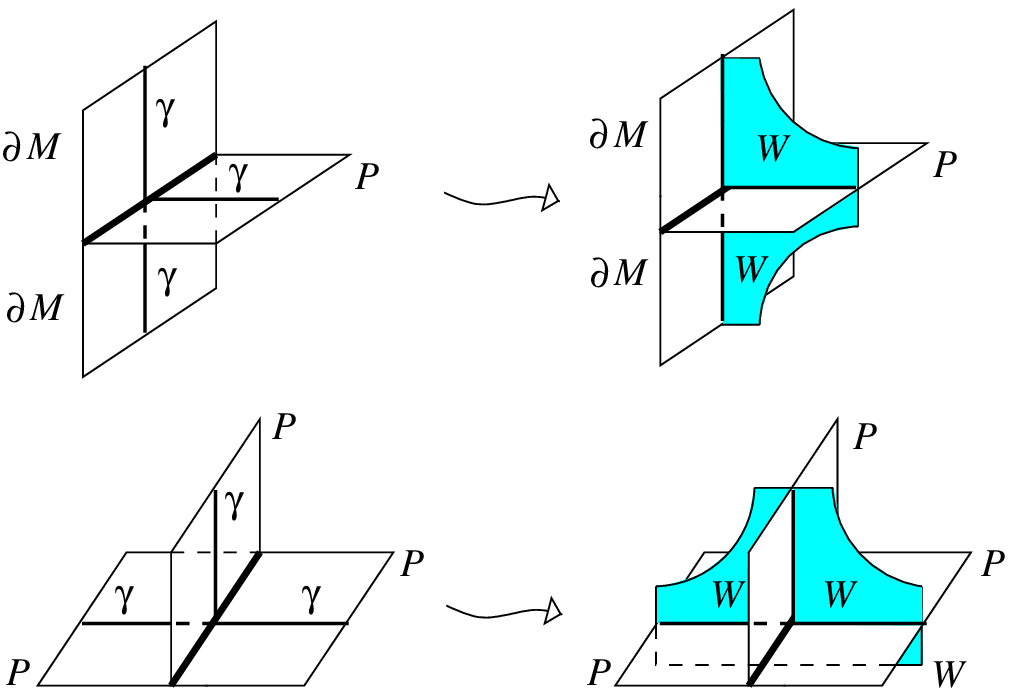}
\nota{Surface determined by a trivalent graph.}\label{traces:fig}
\end{center}\end{figure}

\begin{lemma}\label{split:lem}
Let $P$ be a quasi-standard skeleton of $(M,X)$ and let $\gamma$
be a trivalent graph contained in $(P\cup\partial M)\setminus V(P\cup\partial M)$,
locally embedded as in Fig.~\ref{traces:fig}-left. Then:
\begin{itemize}
\item There exists a properly embedded surface $S$ in $M$ such that
$S\cap (P\cup\partial M)=\gamma$ and $S\setminus \gamma$ is a union of discs. Moreover
$S$ is separating in $M$ if and only if $\gamma$ is separating on $P\cup\partial M$.
\end{itemize}
Assume now that $\gamma\in\{\theta,\sigma\}$ is contained in $P$,
that $S$ is separating in $M$, and that $S\setminus\gamma$ is one disc only. Then:
\begin{itemize}
\item Cutting $M$ along $S$ and choosing $\gamma$ as a spine for the two new
boundary components we get a decomposition $(M,X)=(M_1,X_1)\oplus(M_2,X_2)$
which, at the level of skeleta, corresponds precisely to the splitting of $P$ along $\gamma$.
\end{itemize}
\end{lemma}

\begin{proof}
We first construct a surface with boundary $W$ which meets $P\cup\partial M$
transversely precisely
along $\gamma$, as suggested in Fig.~\ref{traces:fig}-right. Now the portion of $\partial W$ which does not lie on $\partial M$ consists
of a finite number of disjoint circles which can be considered to lie on the boundary
of a concentric sub-ball $B'$ of the ball $M\setminus(P\cup\partial M)$.
These circles bound disjoint discs in $B'$, and if we attach these discs to $W$
we get the desired $S$. Such a $S$ is separating if and only if $\gamma$ is because
any arc in $M\setminus S$ with ends on $P\cup\partial M$
can be homotoped to an arc on $(P\cup\partial M)\setminus\gamma$.
This proves the first assertion. The second assertion is obvious.
\end{proof}

\paragraph{Minimal skeleta are standard}
We will now show a theorem on which most of our results will be based.
We first make an easy remark and then state and prove the theorem,
which implies in particular Proposition~\ref{compl:zero:prop}. Later we will show
Theorem~\ref{naturality:teo}.

\begin{rem}\label{standard:is:proper:rem}
\emph{If $P$ is a nuclear and
standard skeleton of $(M,X)$ then it is properly embedded, namely
$\partial P =\partial M \cap P = X$, and $P\cup\partial M$ is standard
without boundary. Moreover $P\cup\partial M$
is a spine of a manifold bounded by one sphere and some tori and Klein bottles, so
$\chi(P\cup\partial M)=1$. Knowing that $S(P\cup\partial M)$ is 4-valent
and denoting by $f(P)$ the number of faces of $P$, we also see that
$f(P)-v(P)=\#X+1$.}\end{rem}

\begin{teo} \label{standard:teo}
Let $(M,X) \in \calX$ be an irreducible and \ptwoirred\ pair, and let $P$
be a minimal spine of $(M,X)$. Then:
\begin{enumerate}
\item If $c(M,X)>0$ then $P$ is standard;
\item If $c(M,X)=0$ and $X=\emptyset$ then $M\in\{S^3,L_{3,1},\matP^3\}$
and $P$ is not standard;
\item If $c(M,X)=0$ and $X\neq\emptyset$ then $(M,X)$ is one of the $B_i^*$ or $Z_k$,
and $P$ is precisely the skeleton described in Section~\ref{skeleta:section},
so $P$ is standard unless $(M,X)$ is $B_1$ or $B'_1$.
\end{enumerate}
\end{teo}

\begin{proof}
Points (1) and (2), in the closed orientable case, are due
to Matveev~\cite{Matveev:new:book}. Point (3), which requires a rather
careful argument and does not have any closed or even orientable analogue, is new.

We first show that if $P$ is not standard then either
$X=\emptyset$ and $M\in\{S^3,L_{3,1},\matP^3\}$, or
$(M,X)\in\{B_1,B'_1\}$ and $P\in\{P_1,P'_1\}$.
Later we will describe standard skeleta without vertices.

If $P$ reduces to one point of course $M=S^3$.
Let us first assume that $P$ is not purely 2-dimensional, so there is segment $e$
contained in the 1-dimensional part of $P$.
We distinguish two cases depending on whether $e$ lies in $\interior{M}$
or on $\partial M$.

If $e\subset\interior M$, we take a small disc $\Delta$ which intersects $e$ transversely
in one point. As in the proof of Lemma~\ref{split:lem} we attach to
$\partial\Delta$ a disc contained in
the ball $M\setminus(P\cup\partial M)$, getting a
sphere $S\subset M$ intersecting $P$ in one point of $e$. By
irreducibility $S$ bounds a ball $B$, and $P\cap B$ is easily seen to
be a spine of $B$. Nuclearity now implies that $P\cap B$ contains
vertices, so $P\setminus B$ is a skeleton of $(M,X)$ with fewer
vertices than $P$. A contradiction.

If $e\subset\partial M$, let $C$ be the component of
$\partial M$ on which $e$ lies. Since on $C$ there is a circle
which meets $\tau$ transversely in one point of $e$,
looking at the ball
$M\setminus(P\cup\partial M)$ again we see that in $M$ there is a properly
embedded disc $D$ intersecting $P$ in a point of $e$. We have now three cases depending
on the type of the pair $(C,\tau)$.
\begin{itemize}
\item If $(C,\tau)=(T,\theta)$ then $D$ is a compressing disc for $T$, so
by irreducibility $M$ is the solid torus. Knowing that $\partial D$ meets
$P$ only in one point it is now easy to show also that $(M,X)=B_1$ and $P=P_1$.
\item If $(C,\tau)=(K,\theta)$ then $e$ must be contained in the edge $e'''$ of
$\theta$ by Lemma~\ref{nuclear:lem}, and the same reasoning shows that $(M,X)=B'_1$ and $P=P'_1$.
\item If $(C,\tau)=(K,\sigma)$ then $e$ must be contained in the edge
$e'''$ of $\sigma$ by Lemma~\ref{nuclear:lem}. The complement in $K$ of
$\partial D$ is now the union of two M\"obius strips.
If we choose any one of these strips and take its union with $D$,
we get an embedded $\matP^2$ in $M$. Being irreducible and
\ptwoirred, $M$ should then be $\matP^3$, but $\partial M\ne\emptyset$:
a contradiction.
\end{itemize}

We are left to deal with the case where $P$ is purely two-dimensional,
so it is quasi-standard, but it is not standard. Let us first suppose that some
2-component $F$ of $P$ is not a disc.
Then either $F$ is a sphere, so
$P$ also reduces to a sphere only, which is clearly impossible
because $M$ would be $S^2\times[0,1]$, or there exists
a loop $\gamma$ in $F$ such that one of the following holds:
\begin{enumerate}
\item $\gamma$ is orientation-reversing on $F$;
\item $\gamma$ separates $F$ in two components none of which is a disc.
\end{enumerate}

We consider now the closed surface $S$ determined by $\gamma$ as in
Lemma~\ref{split:lem}, and note that $S$ is either $S^2$ or $\matP^2$.
If $S=\matP^2$ we deduce that $(M,X)=\closedpair{\matP^3}$.
If $S=S^2$ irreducibility implies that $S$ bounds a
ball $B$ in $M$. This is clearly impossible in
case (1), so we are in case (2). Now we note that $P\cap B$ must be a nuclear spine of $B$.
Knowing that $F\cap B$ is not a disc it is easy to deduce that $P\cap B$ must contain
vertices. This contradicts minimality because we could replace the whole of
$P\cap B$ by one disc only, getting another skeleton of $(M,X)$ with fewer vertices.

If $P$ is quasi-standard and its 2-components are discs then either $P$ is standard
or $S(P)$ reduces to a single circle. Then it is easy to show that
$P$ must be the triple hat and $(M,X)=\closedpair{L_{3,1}}$.

We are left to analyze the case where $P$
is standard and $c(M,X)=0$, so $X\ne\emptyset$.
Denoting $\#X$ by $n$, Remark~\ref{standard:is:proper:rem} shows that $P$ has $n+1$ faces.

We consider first the case $n=1$. Since $P$ has one edge and two faces, it is easy to
see that it must be homeomorphic to either
$P_2$ or $P_2'$ (see Fig.~\ref{p2spines:fig}) as an abstract polyhedron.
This does not quite imply that $(M,X)$ is $B_2$ or $B'_2$, because in general
a skeleton $P$ alone is not enough to determine a pair $(M,X)$.
However $P\cup\partial M$ certainly does determine $(M,X)$, because it is a standard
spine of $M$ minus a ball, and $X=P\cap\partial M$. We are
left to analyze all the polyhedra of the form $P_2\cup_\psi T$ for
$\psi:\partial P_2\to\theta\subset T$, of the form
$P_2\cup_\psi K$ for $\psi:\partial P_2\to\theta\subset K$,
and of the form $P'_2\cup_\psi K$ for
$\psi:\partial P'_2\to\sigma\subset K$.
Among these polyhedra we must
select those which can be thickened to manifolds with two boundary components
(a sphere plus either a torus or a Klein bottle).
The symmetries of $(T,\theta)$, $(K,\theta)$, and
$(K,\sigma)$ described in Propositions~\ref{torus:spines:prop}
and~\ref{Klein:spines:prop} imply that there are
actually not many such polyhedra.
More precisely, there is just one
$P_2\cup_\psi T$, which gives $B_2$.
There are two $P_2\cup_\psi T$, one of them is not
thickenable (\emph{i.e.}~it is not the spine of any  manifold),
and the other one can
be thickened to a manifold with three boundary components
(a sphere and two Klein bottles).
Finally, there are two $P'_2\cup_\psi K$,
one is not thickenable and the other one gives $B'_2$.
This concludes the proof for $n=1$.

Having worked out the case $n=1$, we turn to $n\ge 2$, so $P$ has $n$ edges and
$n+1$ faces. If a face of $P$ meets $\partial M$ in one arc only, then it meets
$S(P)$ in one edge only, and this edge joins a component of $\partial M$ to
itself, which easily implies that $n=1$, against the current assumption. If a
face of $P$ is an embedded rectangle, with two opposite edges on $\partial M$
and two in $S(P)$, then it readily follows that $n=2$ and $P$ is either
$\theta\times[0,1]$ or $\sigma\times[0,1]$.
As above, to conclude that
$(M,X)\in\{B_0,B'_0,B''_0\}$, we must
consider the various polyhedra obtained by attaching $(T,\theta)$,
$(K,\theta)$, and $(K,\sigma)$ to the upper and lower bases of
$\theta\times[0,1]$ and $\sigma\times[0,1]$. Using again
Propositions~\ref{torus:spines:prop} and~\ref{Klein:spines:prop}
one sees that there are only six such polyhedra.
Three of them are not thickenable, and the other three give
$B_0,B'_0,B''_0$.

Back to the general case with $n\ge 2$, we note that there is a total of $3n$
edges on $\partial M$, so there are $3n$ germs of faces starting from $\partial M$.
Knowing that there is a total of $n+1$ faces and none of them uses one germ only,
we see that at least one face uses two germs only, so it is a rectangle $R$,
possibly an immersed one. If $n=2$ we have three rectangles, one of which must be embedded,
and we are led back to a case already discussed. If $n\ge 3$ then $R$
must be immersed, so in particular it joins a component $(K_1,\sigma_1)$ of
$(\partial M,X)$ to another $(K_2,\sigma_2)$ component. A regular neighborhood in $P$
of $R\cup\sigma_1\cup\sigma_2$ is shown in Fig.~\ref{pkcut:fig}.
\begin{figure}\begin{center}
\figfatta{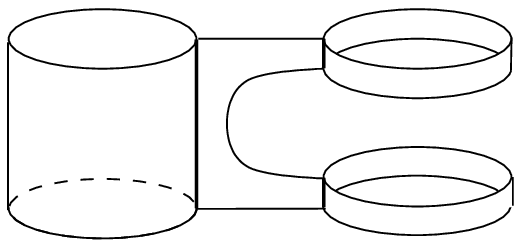}
\nota{An immersed rectangle joins two $(K,\sigma)$ components.}\label{pkcut:fig}
\end{center}\end{figure}
The boundary of this neighborhood is again a graph $\sigma$
which determines a separating Klein bottle according to Lemma~\ref{split:lem}.
If we cut $P$ along $\sigma$ we get a disjoint union
$P''_2\sqcup P'$, which at the level of manifolds gives a splitting
$(M,X)=B''_2\oplus(M',X')$. Moreover $P'$ is a nuclear skeleton of $(M',X')$,
so $c(M',X')=0$, $P'$ is minimal, and $\#X'=n-1$. Now either $(M',X')=B''_0$ and
$P'=P''_0$ or we can proceed,
eventually getting that $(M,X)=B''_2\oplus\ldots\oplus B''_2$,
so $(M,X)=Z_k$ for some $k\ge 3$, and $P$ is the corresponding
skeleton constructed in Section~\ref{skeleta:section}.
The proof is now complete.
\end{proof}

\dimo{Theorem~\ref{naturality:teo}}
By the previous result, a minimal spine of $\closedpair{M}$ is standard with vertices,
and dual to it there is a singular triangulation with $c(M)$
tetrahedra (and one vertex). A singular triangulation of $M$ with $n$ tetrahedra
and $k$ vertices dually gives a standard polyhedron $Q$ with $n$ vertices
whose complement is a union of $k$ balls.
If we puncture $k-1$ suitably chosen faces
of $Q$ we get a skeleton
of $(M,\emptyset)$, whence the conclusion at once.
\finedimo

\section{Finiteness}\label{finiteness:section}

The proof of Theorem~\ref{finiteness:teo}
will be based on the following result.

\begin{prop} \label{standard:SP:connected:prop}
Let $(M,X)$ be an irreducible and \ptwoirred\ pair
such that $c(M,X)>0$ and $(M,X)$ does not split as an assembling
$(M,X)=(N,Y)\oplus B''_2$.
Let $P$ be a standard skeleton of $(M,X)$. Then
every edge of $P$ is incident to at least one vertex of $P$.
\end{prop}

\begin{proof}
Assume by contradiction
that an edge $e$ of $P$ is not incident to any vertex of $P$, \emph{i.e.}~that
both the ends of $e$ lie on $\partial M$. If the ends of $e$
lie on the same spine $\tau\in X$ then
$\tau\cup e$ is a connected component of $S(P)\cup\partial M$.
Standardness of $P$ implies that $P$ has no vertices, which contradicts the
assumption that $c(M,X)>0$. So the ends of $e$ lie on distinct
spines $\tau,\tau'\in X$.
Let $C$ and $C'$ be the components of $\partial M$ on which $\tau$ and $\tau'$ lie,
and let $R$ be a regular neighborhood in $P$ of $C\cup C' \cup e$.
By construction $R$ is a quasi-standard polyhedron with boundary
$\partial R=\tau\sqcup\tau'\sqcup\gamma$. Here $\gamma$ is a trivalent graph
with one component homeomorphic to $\theta$ or to $\sigma$,
and possibly another component homeomorphic to the circle.

Let us first consider the case where $\gamma$ has a circle component $\gamma_0$.
This circle lies on $P$ and is disjoint from $S(P)$. Standardness of $P$
then implies that $\gamma_0$ bounds a disc $D$ contained in $P$ and disjoint
from $S(P)$. In this case we set $\gamma'=\gamma\setminus\gamma_0$ and $R'=R\cup D$.
In case $\gamma$ is connected we just set $\gamma'=\gamma$ and $R'=R$.
In both cases we have found a graph $\gamma'$ homeomorphic to
$\theta$ or to $\sigma$ which separates $P$. Moreover one component
$R'$ of $P\setminus\gamma'$ is standard without vertices and is bounded by
$\tau\sqcup\tau'\sqcup\gamma'$.

According to Lemma~\ref{split:lem}, the graph $\gamma'$
determines a separating surface $S$ in $M$ such that $S\cap P=\gamma'$.
Since $\chi(\gamma')=-1$ and $S\setminus\gamma'$ consists of discs, we have
$\chi(S)\ge 0$. Of course $\chi(S)\ne 1$, for otherwise $S$ would be an embedded
$\matP^2$, but we are assuming that $M$ is irreducible and \ptwoirred\ and has
non-empty boundary. We will now show that if $\chi(S)=2$ then $c(M,X)=0$,
and if $\chi(S)=0$ then $(M,X)$ splits as $(M,X)=(N,Y)\oplus B''_2$.
This will imply the conclusion.

Assume that $\chi(S)=2$, so $S$ is a sphere. We denote by $B$ the open 3-ball
$M\setminus(P\cup\partial M)$ and note that $S\cap B=S\setminus\gamma'$
consists of three disjoint open 2-discs, which cut $B$ into four open 3-balls.
By irreducibility, $S$ bounds a closed 3-ball $D$, and $B\setminus D$
is the union of some of the four open 3-balls just described.
Viewing $(D,\gamma')$ abstractly we can now easily construct a new simple
polyhedron $Q\subset D$ without vertices such that $Q\cap S=\gamma'$
and $D\setminus Q$ consists of three distinct 3-balls, each incident
to one of the three open 2-discs which constitute $S\setminus\gamma'$.
Let us consider now the simple polyhedron $P'=R'\cup_{\gamma'} Q$
viewed as a subset of $M$.
By construction $P'\cap\partial M=\tau\cup\tau'=X$. Moreover
$M\setminus(P'\cup\partial M)$ is obtained from $B\setminus D$
(which consists of open 3-balls) by attaching each of the three 3-balls
of $D\setminus Q$ along only one 2-disc (a component of $S\setminus\gamma'$).
It follows that $M\setminus(P'\cup\partial M)$ still consists
of open 3-balls. By puncturing some of the 2-components
of $P'$ we can then construct a skeleton of $(M,X)$
without vertices, so indeed $c(M,X)=0$.

Assume now that $\chi(S)=0$, so $S$ is a separating torus or Klein bottle.
Lemma~\ref{split:lem} now shows that $(M,X)$ is obtained by assembling
some pair $(N,Y)$ with a pair $(N',Y')$ which has skeleton $R'$.
By construction $R'$ is standard without vertices
and $\partial N'$ has three components, and it was shown within
the proof of Theorem~\ref{standard:teo} that $(N',Y')$ must then
be $B_2''$. This completes the proof.
\end{proof}

\begin{cor} \label{limitazione:su:X:cor}
Let $(M,X)$ be irreducible and \ptwoirred. Assume
$c(M,X)>0$ and there is no splitting
$(M,X)=(N,Y)\oplus B''_2$. Then $\#X\le 2c(M,X)$.
\end{cor}
\begin{proof}
A minimal skeleton $P$ of $(M,X)$ is standard by
Theorem~\ref{standard:teo}, and we have just
shown that each edge of $P$ joins either $V(P)$ to itself or $V(P)$
to $X$. Since $P$ has $c(M,X)$ quadrivalent vertices,
there can be at most $4c(M,X)$ edges reaching $X$.
Each component of $X$ is reached by precisely two edges,
so there are at most $2c(M,X)$ components.
\end{proof}

\dimo{Theorem~\ref{finiteness:teo}} The result is valid for $n=0$ by the
classification carried out in Theorem~\ref{standard:teo}, so we assume $n>0$.
Let $\calF_n$ be the set of all irreducible and \ptwoirred\ pairs $(M,X)$ which
cannot be split as $(M',X')\oplus B''_2$. By Theorem~\ref{standard:teo}, each
such $(M,X)$ has a minimal standard spine $P$ with $n$ vertices.   By
Corollary~\ref{limitazione:su:X:cor}, we have that  $S(P\cup\partial M)$ is a
quadrivalent graph with at most $3n$ vertices. Since $P\cup\partial M$ is a
standard polyhedron, there are only finitely many possibilities for
$P\cup\partial M$ and hence for $(M,X)$.

Given an irreducible and \ptwoirred\ pair $(M,X)$ with
$c(M,X)=n$, either $(M,X)\in\calF_n$ or $(M,X)$ splits along a Klein bottle
$K$ as $(M',X')\oplus B''_2$. The only case where $K$ is compressible
in $M$ is when $(M',X')=B'_2$, but
$B'_2\oplus B''_2=B''_0$ and $c(B''_0)=0$. So $K$ is incompressible,
whence $M'$ is irreducible and \ptwoirred. Moreover $c(M',X')=n$ by
Remark~\ref{B2sec:is:sharp:rem}
(which depends on the now proved Propositions~\ref{compl:zero:prop}
and~\ref{subadditivity:prop}).
Since $(M',X')$ has one boundary component less than $(M,X)$,
we can iterate the process of splitting copies of $B_2''$ only
a finite number of times, and then we get to an element of $\calF_n$.
\finedimo

\section{Additivity}\label{additivity:section}
In this section we prove additivity under connected sum. This will require the theory of
normal surfaces and more technical results on skeleta. We start with an easy general
fact on properly embedded polyhedra.

\begin{prop} \label{separ:surf:prop} Given a pair $(M,X)\in\calX$, let
$Q\subset M$ be a quasi-standard polyhedron with $Q\cap\partial M=\partial
Q\subset X$. Assume that $M\setminus Q$ has two components $N'$ and $N''$.
Then the $2$-components of $Q$ that separate $N'$ from $N''$ form a closed
surface $\Sigma(Q) \subset Q \subset \interior{M}$ which cuts $M$ into two
components. \end{prop}

\begin{proof} Let $e$ be an edge of $Q$, and let $\{F_1, F_2, F_3\}$ be
the triple of (possibly not distinct) faces of $Q$ incident to $e$. The
number of $F_i$'s that separate $N'$ from $N''$ is even; it follows that
$\Sigma(Q)$ is a surface away from $V(Q)\cup\partial Q$. Let $C$ be a
boundary component of $M$, containing $\tau\in X$. Since
$C\setminus \tau$ is a disc, which is adjacent either to $N'$ or to
$N''$ (say $N'$), then each $2$-component of $Q$ incident to $\tau$ has $N'$ on both
sides. So $\Sigma(Q)$ is not adjacent to $\partial Q$. Finally, since
$\Sigma(Q)$ intersects the link of each vertex either nowhere or in a
loop, then $\Sigma(Q)$ is a closed surface.
It cuts $M$ in two components because $N'$ and $N''$ lie on opposite sides of
$\Sigma(Q)$. \end{proof}

\paragraph{Normal surfaces} Given a pair $(M,X)\in\calX$, let
$P$ be a nuclear skeleton of $(M,X)$. The simple polyhedron
$P\cup\partial M$ is now a spine of $M$ with a ball $B\subset M$
removed. Choose a triangulation of $P\cup\partial M$, and let $\xi_P$
be the handle decomposition of $M\setminus B$ obtained by thickening the
triangulation of $P\cup\partial M$, as in~\cite{Matveev:new:book}. In this
paragraph we will study normal spheres in $\xi_P$. Note that
there is an obvious one, namely the sphere parallel to $\partial B$ and
slightly pushed inside $\xi_P$. The following result deals with the other normal spheres.
Its proof displays another remarkable difference between the orientable and
the general case. Namely, it was shown in~\cite{MaPe:EM} that, when $M$ is orientable,
any normal surface reaching $\partial M$ actually contains a component of $\partial M$.
On the contrary, when $(\partial M,X)$ contains some $(K,\sigma)$
component, an arbitrary normal surface can reach $K$ without containing it.
As our proof shows, however, this cannot happen when the surface is a sphere.

\begin{prop} \label{normal:sphere:prop}
Let $P$ be a nuclear skeleton of $(M,X)\in\calX$, and let $S$ be a normal sphere
in $\xi_P$. Then:
\begin{itemize}
\item There exists a simple polyhedron $Q$ such that $v(Q)\le v(P)$,
$Q\cap\partial M=X$ and $M\setminus(Q\cup\partial M)$ is a regular neighborhood
of $S$.
\end{itemize}
Suppose now in addition that $P$ is standard, that $c(M,X)>0$ and
that $S$ is not the obvious sphere $\partial N(P\cup\partial M)$. Then:
\begin{itemize}
\item There exists $Q$ as above with $v(Q)<v(P)$.
\end{itemize}
\end{prop}

\begin{proof} Every region $R$ of $P$ carries a color $n\in\matN$ given by
the number of sheets of the local projection of $S$ to $R$. Now we cut
$P\cup\partial M$ open along $S$ as explained in~\cite{Matveev:new:book},
\emph{i.e.}~we replace each $R$ by its $(n+1)$-sheeted cover contained in the
normal bundle of $R$ in $M$. As a result we get a polyhedron $P'\subset M$
which contains $\partial M$, such that $M\setminus P'$ is the disjoint
union of an open ball $B$ and an open regular neighborhood $N$ of $S$ in
$M$. By removing from each boundary component $C\subset \partial M$ the
open disc $C\setminus\tau$ we get a polyhedron $P''$ intersecting
$\partial M$ in $X$. Now we puncture a $2$-component which separates $B$
from $N$ and claim that the resulting polyhedron $Q$ is as
desired. Only the inequalities between $v(P)$ and $v(Q)$ are non-obvious.

We first prove that all the vertices of $P\cup\partial M$ which lie on $\partial M$
disappear either when we cut $P$ along $S$ getting $P'$ or later when we remove
$\partial M\setminus X$ from $P'$ to get $P''$. This of course implies the first
assertion of the statement. We concentrate on one component $(C,\tau)$ of $(\partial M,X)$.
By Lemma~\ref{nuclear:lem} either both vertices of $\tau$ are vertices of
$P\cup\partial M$ or none of them is. In the latter case there is nothing to show,
so we assume that there are three (possibly non-distinct) 2-components
of $P$ incident to $\tau$. Let $v$ and $v'$ be the vertices of $\tau$.
Looking first at $v$, we denote by $(n,n,n,p,q,r)$ the colors
of the six germs at $v$ of 2-component of $P\cup\partial M$. Here $n$ corresponds
to $C\setminus\tau$, which is triply incident to $v$.

The compatibility equations of normal surfaces now readily imply that
that (up to permutation) $r$ is even, $p=q\ge r$, and that $n\ge p/2$ when $p=q=r$. Moreover:
\begin{itemize}
\item $v$ disappears in $P'$ if $p=q>r$;
\item $v$ survives in $P'$ and remains on $\partial M$,
so it disappears in $P''$, if $p=q=r$ and $n=p/2$;
\item $v$ survives in $P'$ and moves to $\interior{M}$ if $p=q=r$ and $n>p/2$.
\end{itemize}
Now if $\tau=\theta$ then the same coefficients appear at $v'$. The only case where
$v$ and $v'$ do not both disappear in $P''$ is when $p=q=r$ and $n>p/2$. But in this case
$S$ would contain $n-p/2$ parallel copies of $C$, which is impossible.
The case $\tau=\sigma$ is easier, because if $v$ survives in $P''$ the situation
is as in Fig.~\ref{mobsphere:fig}.
\begin{figure}\begin{center}
\figfatta{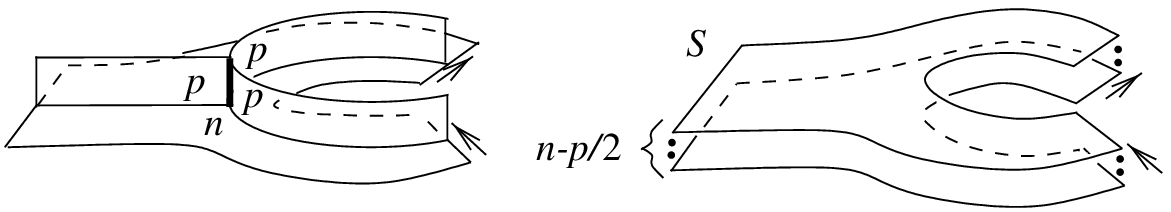}
\nota{M\"obius strips in a normal surface.}\label{mobsphere:fig}
\end{center}\end{figure}
This is absurd because $S$ would contain M\"obius strips.

Now we turn to the second assertion. If $v(P'')<v(P)$ the conclusion is
obvious, so we proceed assuming $v(P)=v(P'')$. It is now sufficient to
show that some face of $P''$ which separates $B$ from $N$ contains
vertices of $P''$, because we can then puncture such a face and collapse
the resulting polyhedron until it becomes nuclear, getting fewer vertices.
Assume by contradiction that there is no such face.

We note that $P''$ is the union of a quasi-standard polyhedron $P'''$ and
some arcs in $X$. The 2-components of $P''$ which separate $B$ from $N$
are the same as those of $P'''$, so they give a closed surface
$\Sigma\subset P''$ by Proposition~\ref{separ:surf:prop}. From the fact
that $v(P'')=v(P)$ we deduce that near a vertex of $P$ the transformation
of $P$ into $P''$ can be described as in Fig.~\ref{localtrasf:fig},
\begin{figure}\begin{center}
\figfatta{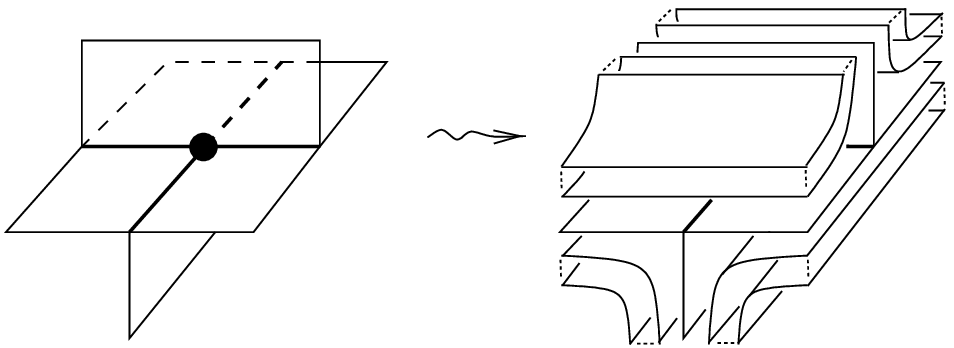}
\nota{Transformation
of $P$ into $P''$ near a vertex of $P$.}\label{localtrasf:fig}
\end{center}\end{figure}
namely $P''$ can be identified near the vertex
with $P\cup S$. Of course this does not imply that globally $P''=P\cup S$,
because the components of $P''$ playing the role of $P$ near vertices may
not match across faces.

The closed surface
$\Sigma$ cannot be disjoint from $S(P'')$, because otherwise $S$ would be
the obvious sphere $\partial B$. On the other hand we are supposing
$\Sigma\cap V(P'')=\emptyset$, so $\Sigma\cap S(P'')$ must be a
non-empty union of loops. In particular, $S(P'')$ contains a loop $\gamma$
disjoint from $V(P'')$.

Figure~\ref{localtrasf:fig}
now shows that $S(P'')$  coincides with $S(P)$ away from $\partial M$. Using
the analysis of the transition from $P$ to $P''$ near $\partial M$
already carried out above, we also see that near a component
$(C,\tau)$ of $(\partial M,X)$ either $S(P'')$ coincides with $S(P)$ or it
is obtained from $S(P)$ by adding one edge of $\tau$,
and then slightly pushing the result inside $M$. When
$(C,\tau)=(K,\sigma)$ the edge added is necessarily $e'''$.
This implies that the loop $\gamma$ described above can be viewed as a loop in
$S(P\cup\partial M)$ such that $\gamma\cap V(P)=\emptyset$. In addition,
if $\gamma$ contains a vertex of $P\cup\partial M$ on a certain
component of $\partial M$ then it contains also the other vertex in that component.
This readily implies that the union of $\gamma$ with all the $\tau$'s in $X$ touched
by $\gamma$ is a connected component of $S(P\cup\partial M)$. But
$P\cup\partial M$ is standard, so $S(P\cup\partial M)$ is connected, and we deduce that
$P$ has no vertices. A contradiction.
\end{proof}

\dimo{Theorem~\ref{additivity:teo}}
We have already noticed that $c(S^2\timtil S^1)=c(S^2\timtil S^1)=0$
and that $c$ is subadditive. Let us consider now a non-prime pair $(M,X)$
and a minimal skeleton $P$ of $(M,X)$. Since $(M,X)$ is not prime, there exists a
normal sphere $S$ in $\xi_P$ which is essential in $M$, namely either it is
non-separating or it separates $M$ into two manifolds both different from $B^3$.
Then we apply the first point of Proposition~\ref{normal:sphere:prop} to $P$ and $S$,
getting a polyhedron $Q$.

If $S$ is separating and splits $(M,X)$ as $(M_1,X_1)\#(M_2,X_2)$,
we must have that $Q$ is the disjoint
union of polyhedra $Q_1$ and $Q_2$, where $Q_i$ is a skeleton of $(M_i,X_i)$.
Since $v(Q_1)+v(Q_2)=v(Q)\le v(P)$ we deduce that
$c(M,X)\ge c(M_1,X_1)+c(M_2,X_2)$, so equality actually holds.

If $S$ is not separating we identify a regular neighborhood of $S$ in $M$
with $S\times (-1,1)$ and note that there must exist a face of $Q$ having
$S\times(-1,-1+\varepsilon)$ on one side and $S\times(1-\varepsilon,1)$ on the
other side. We puncture this face getting a polyhedron $Q'$. Now $Q'$ is a
skeleton of a pair $(M',X)$ such that $(M,X)=(M',X)\# E$ where $E$ is
$S^2\times S^1$ or $S^2\timtil S^1$. Moreover $v(Q')=v(Q)\le v(P)$,
hence $c(M,X)\ge c(M',X)$, so equality actually holds.

We have shown so far that an essential normal sphere in $(M,X)$ leads to a
non-trivial decomposition $(M,X)=(M_1,X_1)\#(M_2,X_2)$ on which complexity
is additive. If $(M_1,X_1)$ and $(M_2,X_2)$ are prime we stop, otherwise
we iterate the procedure until we find one decomposition of $(M,X)$ into
primes on which complexity is additive. Since any other decomposition into
primes actually consists of the same summands, we deduce that complexity
is always additive on decompositions into primes. If we take the connected
sum of two non-prime manifolds then a prime decomposition of the result is
obtained from prime decompositions of the summands, so additivity holds
also in general.
\finedimo

\section{Sharp assemblings}\label{assemblings:section}

In this section we prove Theorem~\ref{sharp:assemble:teo}.

\paragraph{Pairs with standard minimal skeleta} The main ingredient for
Theorem~\ref{sharp:assemble:teo} is the following partial converse of
Theorem~\ref{standard:teo}:

\begin{teo}\label{inverse:standard:teo}
If a pair $(M,X)$ has a standard minimal skeleton then it is irreducible.
\end{teo}

\begin{proof} If $c(M,X)=0$ the conclusion follows from the classification
of standard skeleta without vertices, which was carried out within the
proof of Theorem~\ref{standard:teo}. So we assume $c(M,X)>0$. We proceed
by contradiction and assume that there exists an essential sphere, whence
a normal one $S$ with respect to a standard minimal skeleton $P$. We can
now apply the second point of Proposition~\ref{normal:sphere:prop} to $P$
and $S$, getting a polyhedron $Q$. By adding an arc to $Q$ we get a new
skeleton of $(M,X)$ with fewer vertices than $P$: a contradiction.
\end{proof}

\paragraph{Exceptional bricks}
We show in this paragraph that the bricks $B_1$ and $B'_1$, which we regard to
be exceptional by Theorem~\ref{standard:teo}, never appear in the splitting
of a positive-complexity irreducible
and \ptwoirred\ pair. This fact will be used in the proof
of Theorem~\ref{sharp:assemble:teo}.

\begin{lemma}\label{b1:assembling:lem} Let $B_1^*\oplus (M,X)=(N,Y)$ be a
sharp assembling with $(M,X)$ irreducible and \ptwoirred. Then $c(M,X)=0$
and $$(N,Y)\in\{S^3,L_{3,1},\matP^3,S^2\times S^1,S^2\timtil S^1\}.$$
\end{lemma}

\begin{proof}
We first assume that $(M,X)$ cannot be expressed as $(M',X')\oplus B''_2$,
we choose a minimal skeleton $P$ of $(M,X)$,
and we apply Propositions~\ref{compl:zero:prop}
and~\ref{standard:SP:connected:prop}, which easily imply that either
$(M,X)=B_i^*$ with $i\le2$ or every face of $P$ contains vertices.
If we attach $P_1^*$ and $P$ along the map which gives the
assembling we get a skeleton $Q$ of $(N,Y)$ having $c(N,Y)$ vertices.
Recall now that $P_1^*$ has a 1-dimensional portion, namely a free segment
$e$ on $\partial M$. If $P$ has vertices we readily deduce that
$Q$ can be collapsed to a subpolyhedron with fewer vertices: a contradiction.
So $(M,X)$ must be of type $B_i^*$ with $i\le2$.
The non-trivial assemblings
$B_1^*\oplus B_i^*$ are easily discussed
and the conclusion follows.

Assume now that $(M,X)=(M',X')\oplus B''_2$. Noting that $B_1^*$ has a
$\theta$ on its boundary, we deduce that $B''_2$ is assembled to $(M',X')$.
Iterating the splitting of copies of $B''_2$ and
applying Remark~\ref{B2sec:is:sharp:rem} and Theorem~\ref{sharp:split:teo}
we get that $(N,Y)=(B_1^*\oplus (M'',X''))\oplus B''_2\oplus\ldots\oplus B''_2$,
where $(M'',X'')$ is irreducible and \ptwoirred\ and cannot
be split as $(M''',X''')\oplus B''_2$, and the assembling
$B_1^*\oplus (M'',X'')$ is sharp.
So $B_1^*\oplus (M'',X'')\in
\{S^3,L_{3,1},\matP^3,S^2\times S^1,S^2\timtil S^1\}$, but
no $B''_2$ can be assembled to any of these manifolds.
\end{proof}

\begin{rem}\emph{By Theorem~\ref{sharp:split:teo}, if we know that
the \emph{result} $(N,Y)$ of a sharp assembling
$B_1^*\oplus(M,X)$ is irreducible and \ptwoirred,
we can apply the previous lemma to deduce that
$(N,Y)\in\{S^3,L_{3,1},\matP^3\}$.} \end{rem}

\paragraph{Faces incident to a spine}
For the proof of Theorem~\ref{sharp:assemble:teo} we need
another preliminary result.

\begin{prop}\label{three:faces:on:tau:prop}
Let $P$ be a standard skeleton of an irreducible and \ptwoirred\ pair $(M,X)$.
Assume that $(M,X)\not\in\{B_1,B'_1,B_2,B'_2\}$. Then for every $\tau\in X$
there are three pairwise distinct faces of $P$ incident to $\tau$.
\end{prop}

\begin{proof}
Let $F$ be doubly incident to
$\tau\subset C\subset\partial M$,
and let $\alpha$ be an arc properly embedded in $F$ with endpoints
on different edges of $\tau$. If we cut
$C$ open along $\tau$ we get a hexagon $H$ as in
Fig.~\ref{hexa:fig}, with identifications which allow to reconstruct $C$.
\begin{figure}\begin{center}
\figfatta{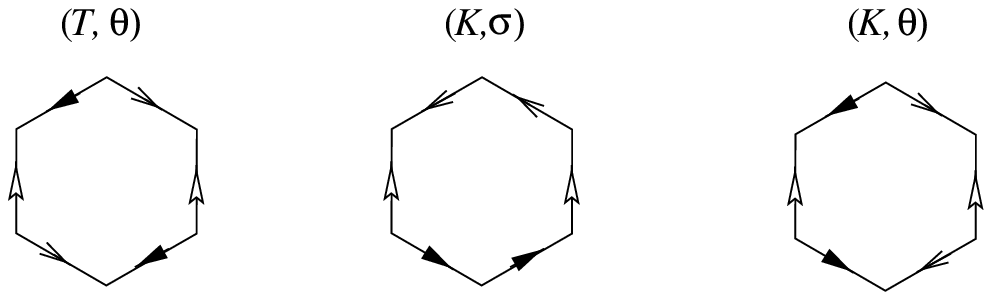}
\nota{Hexagons.}\label{hexa:fig}
\end{center}\end{figure}

The two endpoints of $\alpha$ give rise on $\partial H$
to four points identified in pairs.
Now we choose along $\alpha$ a vector field transversal
to $F$, and we examine this vector at the four points on $\partial H$.
At two of the four points the vector
will be directed towards the interior of $H$, and we join these two
points by an arc $\beta_1$ properly embedded in $H$.
We also join the other two points by
another arc $\beta_2$ and arrange that
$\beta_1$ and $\beta_2$ intersect transversely in at most one
point. Now $\alpha\cup\beta_i$ is a loop for $i=1,2$ and,
as in the proof of Lemma~\ref{split:lem}, we see that
$\alpha \cup \beta_i$ bounds a disc $D_i$ in $M$.
This easily implies that $\beta_1\cap\beta_2$ actually must be empty, for
otherwise $D_1$ and $D_2$ would give rise, in the complement $B^3$
of a regular neighborhood $P\cup\partial M$, to two proper
discs whose boundaries
intersect only once and transversely.

Since $\beta_1\cap\beta_2$ is empty, $D_1\cup D_2$ is a disc properly
embedded in $M$, and the boundary $\beta_1\cup\beta_2$ of this disc
is essential in $C$, because it intersects $\tau$ in two distinct edges.
By irreducibility, $M$
is a solid torus or a solid Klein bottle. If it is a solid torus, since
$\tau=\theta$ meets the meridional disc in two points only, it readily follows
that $(M,X)$ is $B_1$ or $B_2$, against the hypotheses. If it is a solid
Klein bottle, then uniqueness of the embedding of $\theta$ and $\sigma$ in
$K$ implies that $(M,X)$ is $B'_1$ or $B'_2$.
\end{proof}

\dimo{Theorem~\ref{sharp:assemble:teo}}
Both when $(N,Y)=(M,X)\oplus(M',X')$ and when $(N,Y)=\odot(M'',X'')$
we have in $N$ a two-sided torus or Klein bottle $C$ cutting
along which we get a (possibly disconnected) irreducible and
\ptwoirred\ manifold. If $C$ is incompressible in $N$
the desired conclusions follow from routine topological
arguments~\cite{Hempel:book}. The only case where $C$ is compressible
is that of an assembling involving a solid torus or Klein bottle.
So we only have to show irreducibility of $N$
when $(N,Y)=(M,X)\oplus(M',X')$.

Take minimal skeleta $P$ and $P'$ of $(M,X)$ and $(M',X')$. The case where
one of $(M,X)$ or $(M',X')$ is $B_1$ or $B'_1$ was already discussed in
Lemma~\ref{b1:assembling:lem}, so by Theorem~\ref{standard:teo} we have
that $P$ and $P'$ are standard. Let the assembling be performed along
boundary components $(C,\tau)$ and $(C',\tau')$. If the three faces of $P$
incident to $\tau$ are distinct, and similarly for $P'$ and $\tau'$, then
gluing $P$ to $P'$ we get a standard minimal skeleton of $(N,Y)$, so
$(N,Y)$ is irreducible by Theorem~\ref{inverse:standard:teo}. Otherwise,
by Proposition~\ref{three:faces:on:tau:prop}, up to permutation we have
$(M,X)\in\{B_1,B'_1,B_2,B'_2\}$. The case $(M,X)=B_1^*$ was already
discussed. If $(M,X)=B_2^*$ but $(M',X')\not\in\{B_1,B'_1,B_2,B'_2\}$,
from the shape of the skeleton $P_2^*$ (see Fig.~\ref{p2spines:fig}) we
deduce again that $(N,Y)$ has a standard minimal skeleton. If
$(M',X')=B_2^*$ then either $(N,Y)$ is a lens space, so it is irreducible,
or it belongs to $\{S^2\times S^1,S^2\timtil S^1\}$.
\finedimo

\appendix
\section{Some facts about the Klein bottle}\label{klein:section}
In this appendix, following Matveev~\cite{Matveev:new:book},
we classify all simple closed loops
on the Klein bottle $K$ and we deduce
Proposition~\ref{Klein:spines:prop} from this classification.
We also mention two more results on $K$ which
easily follow from the classification.
These results are strictly speaking not necessary
for the present paper, and they are probably well-known to experts,
but we have decided to include them because they
show a striking difference which exists between the orientable
and the non-orientable case.

\begin{prop}\label{Klein:loops:prop}
There exist on the Klein bottle only four non-trivial loops up to isotopy,
as shown in Fig.~\ref{kleinloops:fig}.
\begin{figure}\begin{center}
\figfatta{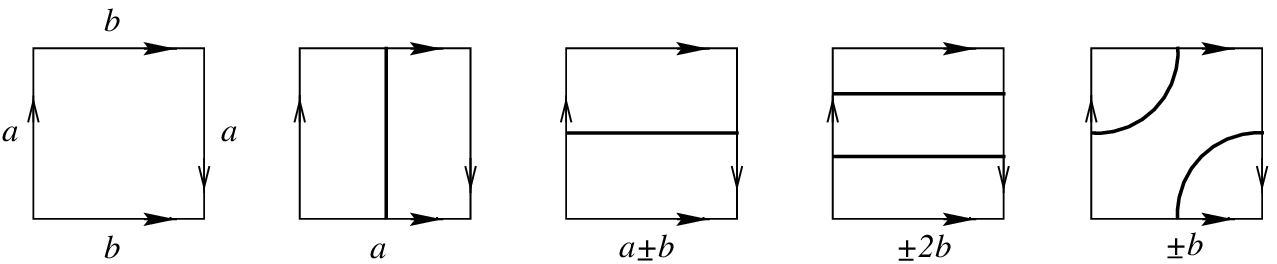}
\nota{Non-trivial loops on the Klein bottle.}\label{kleinloops:fig}
\end{center}\end{figure}
These loops are determined
by their image in $H_1(K;\matZ)=\langle a,b| a+b=b+a,\ 2a=0\rangle$, as
also shown in the picture.
Moreover $a$ and $\pm2b$ are orientation-preserving on $K$, while
$\pm b$ and $a\pm b$ are orientation-reversing.
\end{prop}

\begin{proof}
A non-trivial loop is isotopic to one which is normal with respect
to a triangulation of $K$, \emph{i.e.}~it appears as in Fig.~\ref{klproof:fig}.
\begin{figure}\begin{center}
\figfatta{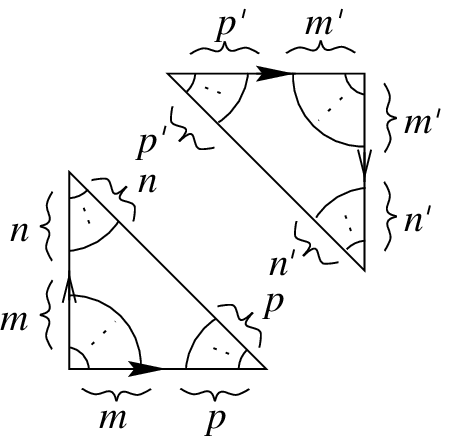}
\nota{Normal loops in a triangulation of $K$.}\label{klproof:fig}
\end{center}\end{figure}
We must have $n+m=n'+m'$, $n+p=n'+p'$, $m+p=m'+p'$, so
$n'=n$, $m'=m$, $p'=p$. If $p>m$, we further distinguish: if
$n<p$, since we look for a connected curve, we get $n=m=0$ and $p=1$, whence
the loop $a$; if $n>p$ we do not get any solution;
if $n=p$ we get $m=0$ and $n=p\in\{1,2\}$, whence the loops $\pm b$ and $\pm 2b$.
If $m>p$ we must have $p=n=0$ and $m\in\{1,2\}$, whence the loops
$\pm b$ and $\pm 2b$ again. If $m=p$,
since the connected curve we look for is also non-trivial,
we must have $m=p=0$ and $n\in\{1,2\}$,
whence the loops $a\pm b$ and $\pm 2b$.
\end{proof}

\dimo{Proposition~\ref{Klein:spines:prop}}
We start by showing that $\sigma$ embeds uniquely as a spine of $K$.
The closed edges $e'$ and $e''$ of $\sigma$ are disjoint simple
loops in $K$, and they must be orientation-reversing.
It easily follows that $\{e',e''\}$ must be
$\{\pm b,a\pm b\}$. Now the ends of $e'''$
can be isotopically slid over $e'$ and $e''$ to reach the
position of Fig.~\ref{twodimspin:fig}-centre, and uniqueness is
proved.

\begin{figure}\begin{center}
\figfatta{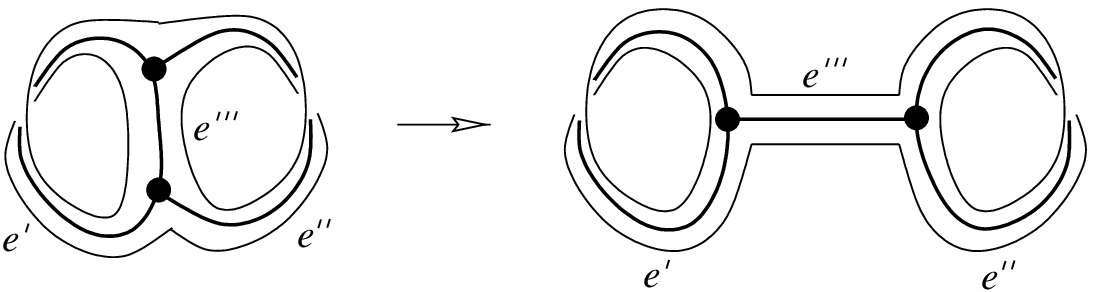}
\nota{A move changing a spine $\theta$ of $K$ into a spine $\sigma$.}\label{spinemove:fig}
\end{center}\end{figure}

Turning to the uniqueness of the embedding of $\theta$, note that two
of the three simple closed loops contained in $\theta$
must be orientation-reversing on $K$. Let $e'''$ be the edge contained
in both these loops. If we perform the move shown in Fig.~\ref{spinemove:fig}
along $e'''$ we get a spine $\sigma$ of $K$, and the newborn
edge is the edge $e'''$ of $\sigma$. So $\theta$ is obtained from
$\sigma$ by the same move along $e'''\subset\sigma$.
The embedding of $\sigma$ being unique, we deduce the same
conclusion for $\theta$.

Having proved uniqueness, we must understand symmetries.
Our description obviously implies that, in both $\sigma$ and $\theta$,
the edges $e'$ and $e''$ play symmetric roles, while the role
of $e'''$ is different, and the conclusion easily follows.
The same conclusion could also be deduced from Fig.~\ref{hexa:fig} or
from Proposition~\ref{Klein:auto:prop} below.
\finedimo

\begin{prop}
If $\solK$ is the solid Klein bottle and $K=\partial\solK$ then
every automorphism of $K$
extends to $\solK$. In particular, there is only
one possible ``Dehn filling''
of a Klein bottle in the boundary of a given manifold.
\end{prop}

\begin{proof}
Proposition~\ref{Klein:loops:prop} shows that
the meridian $a$ of $\solK$ can be characterized in
$K=\partial\solK$ as the only orientation-preserving
loop having connected
complement. So every automorphism of $K$ maps the meridian
to itself and the conclusion follows.
\end{proof}

\begin{prop}\label{Klein:auto:prop}
The mapping class group of $K$ is isomorphic to
$\matZ/_{2\matZ}\times\matZ/_{2\matZ}$ and
every automorphism of $K$ is determined up to
isotopy by its action on $H_1(K;\matZ)$.
\end{prop}

\begin{proof}
It is quite easy to construct commuting
order-2 automorphisms $\phi$ and $\psi$ of $K$ such that
their action on $H_1(K;\matZ)$ is given by
$$\phi(a)=a,\quad\phi(b)=-b,\qquad\psi(a)=a,\quad\psi(b)=a+b.$$
Given any other automorphism $f$, combining the
geometric characterization of $a$ with
the observation that $a$ is isotopic (not only homologous)
to itself with
opposite orientation, we deduce that (up to isotopy) $f$ is the identity
on $a$. Up to composing $f$ with $\phi$ we can assume that
$f$ is actually the identity also near $a$, so $f$ restricts
to an automorphism of the annulus $K\setminus a$ which is the
identity on the boundary. The mapping class group relative
to the boundary of the annulus is now infinite cyclic generated by
the restriction of $\psi$ (but $\psi$ has order 2 when
viewed on $K$), and the conclusion follows.
\end{proof}

\end{document}